\documentclass[preprint,12pt]{elsarticle}

\usepackage{amssymb}
\usepackage{amsmath}
\usepackage{amsthm}
\usepackage{amsfonts}
\usepackage{mathrsfs}
\usepackage{bm}
\usepackage{setspace}
\usepackage{graphicx}
\usepackage{caption}
\usepackage{subcaption}
\usepackage[titletoc,toc,title]{appendix}
\usepackage{rotating}
\usepackage[pdftex,hypertexnames=false,linktocpage=true]{hyperref}
\hypersetup{colorlinks=true,linkcolor=blue,anchorcolor=blue,citecolor=blue,filecolor=blue,urlcolor=blue,bookmarksnumbered=true,pdfview=FitB}

\journal{Renewable Energy}

\begin{document}

\begin{frontmatter}



\title{Constructing low-dimensional stochastic wind models through hierarchical spatial temporal decomposition}


\author[ME]{Qiang Guo\fnref{fn1}}
\author[GAS]{D. Rajewski} 
\author[GAS]{E. Takle}
\author[ME]{Baskar Ganapathysubramanian\corref{cor}}

\address[ME]{Department of Mechanical Engineering, Iowa State University, Ames, Iowa}
\address[GAS]{Department of Geological and Atmospheric Sciences, Iowa State University, Ames, Iowa}
\cortext[cor]{Corresponding author: Tel.:+1 515-294-7442; fax: +1 515-294-3261. E-mail address: baskarg@iastate.edu, URL: http:/www3.me.iastate.edu/bglab/ }
\fntext[fn1]{This author is currently employed at NGC Transmission Equipment (America), Inc.}

\begin{abstract}
Current wind turbine simulations are usually based on simplified wind load models that assume steady wind speed and constant wind speed gradient profile. While various turbulence generating methods have been successfully used for decades, they lack of the ability to reproduce variabilities in wind dynamics, as well as inherent stochastic structures (like temporal and spatial coherences, sporadic bursts, high shear regions). It is important to incorporate these stochastic properties in analyzing and designing the next generation wind-turbines. This necessitates a more realistic parameterization of the wind that encodes location-, topography-, diurnal-, seasonal and stochastic affects. However, such a parameterization is useful in practice only if it is relatively simple (low-dimensional). Interestingly, the data to construct such models are available at various resolutions from meteorological measurements. Such meteorology data have rarely been used in wind turbine design in 
spite of the fact that they contain rich information about location- and topography specific wind speed data. In this work, we develop a hierarchical temporal and spatial decomposition (TSD) of large-scale meteorology data to construct a low-dimensional yet realistic stochastic wind flow model. The TSD framework is based on a two-stage model reduction: Bi-orthogonal Decomposition (BD) followed by Karhunen-Loeve expansion KLE). We showcase this framework on data from a recent meteorological study called CWEX-11 (Crop Wind-energy EXperiment 2011). Starting from this large data set, a low dimensional stochastic parameterization is constructed using the developed framework. The resulting time-resolved stochastic model encodes the statistical behavior exhibited by the actual wind flow. Moreover, the temporal modes encode the variation of wind speed in the mean sense and resolve diurnal variations and temporal correlation while the spatial modes provide deeper insight into spatial coherence of the wind field - which is a key aspect in current wind turbine sizing, design and classification. Comparison of several important turbulent properties between the simulated wind flow and the original dataset show the utility of the framework. We envision this framework as a useful complement to existing wind simulation codes.
\end{abstract}


\end{frontmatter}


\section{Introduction}\label{intro}
The reliability analysis of wind turbines requires a rigorous incorporation of the effect of randomness in the wind load. The wind load on the turbine is a stochastic process whose direction and speed depends on location, and time. In other words, `` ... at any instant they are distributed irregularly in space, at any point in space they fluctuate chaotically in time, and at given point and a given time they vary randomly from realization to realization."~\cite{Wyngaard2010} This randomness in wind conditions leads directly to the fluctuations in rotating speed of the rotor (for variable-speed wind turbine). Furthermore, randomness of wind direction causes random yaw motion of nacelle, which can further affect wind turbine's productivity. Turbulence introduces unsteady loads on the blades. All of these effects excite structural vibrations, introduce components failures, and further reduce the lifespan of wind turbine. Thus, stochastic wind loads have a critical impact on wind turbine performance. To improve 
reliability, a crucial step is to have a better understanding of how they are affected by stochastic wind loads. This can be accomplished using stochastic analysis. Stochastic analysis provides strategies to evaluate the relationship between a stochastic input and various quantities of interest. In general, stochastic analysis is done by modeling the system of interest, incorporating uncertainties in system model or its inputs, and getting probability distributions of the system output. For instance, to quantify the effect of the uncertainty that is introduced by the random wind conditions, stochastic analysis (like polynomial chaos~\cite{Xiu2002,Xiu2003}, stochastic collocation~\cite{Ganapathysubramanian2007}, or Monte Carlo analysis~\cite{Ghanem1998}) on a wind turbine model can be performed~\cite{larsen2007}.  {\it A key ingredient to successful stochastic analysis is to have a realistic, viable models of the input variability} -- in this case, stochastic wind models that encode location-, topography-, diurnal-, and 
seasonal affects.

Most of the past wind turbine design work was based on simplified wind load models that assume steady wind speed, constant wind speed gradient profile and constant turbulence intensity~\cite{Quarton1998}. This kind of simplified wind model does not provide any insight into the effect of randomness in the wind load. They have, however, been used with some success in the context of {\it deterministic wind turbine} analysis~\cite{Bazilevs2011,Bazilevs2011a}. Subsequently, wind models that are capable of describing wind stochasticity have been developed. Gaussian and Weibull distributions have been used in approximating histograms of hourly mean wind speeds. These methods are based on parameter estimation techniques such as maximum likelihood method (MLM) (see Carta et al.~\cite{Carta2009}). While useful is representing randomness, these techniques lack the ability to describe {\it temporal and spatial} correlations of wind, or location-dependent characterization capabilities. This results in the development of 
methods that incorporate spatial variations, spatial correlation, and temporal correlation.

There have been several seminal works that deal with preserving either spatial correlation or temporal correlation. One of such examples is Veers/Sandia method~\cite{veers1988}. The method first uses empirical coherence function and wind speed power spectral densities (PSDs) to describe the wind field. Choleski decomposition is then used to decompose the cross spectrum matrix. Velocity time-series at different locations are finally calculated with a inverse Fourier transform process. Another example is Mann's method~\cite{mann1998}, where the wind field is described in wave vector space and then transformed into spatial domain. Recently, proper orthogonal decomposition (POD) was used~\cite{saranyasoontorn2008} to simulate wind turbine inflow. In this method, stochastic wind flow is viewed as wind field snapshots at different instances. A spatial covariance matrix is decomposed to 
get multiple characteristic modes that are further used in constructing a stochastic wind field. By performing POD analysis at each time step, the wind flow is reproduced. Although these three methods are easy to implement and efficient enough to get random wind flows, they do not account for {\it both temporal correlation and spatial correlation of the wind flow at the same time}. In the software package TurbSim~\cite{jonkman2009turbsim}, a turbulence simulation code developed by National Renewable Energy Laboratory (NREL), Veers method is used. However, due to its limitation in describing temporal coherence, additional information about coherent structure has to be included in the model so that the turbulent structure in atmosphere is more accurately represented~\cite{kelley2005impact}.

The motivation for this work is the fact that no single method exists that can provide a data-driven model that seamlessly accounts for spatial variations, diurnal variations, and temporal correlations \footnote{The methods discussed here are invariably parametric --  where the model (and parameters) are set based on priori assumptions. For example, the probability density function of the velocity fluctuations in homogeneous turbulence is assumed to be Gaussian or sub-Gaussian~\cite{jimenez1998}. Although this assumption is often valid for many physical phenomena (including wind on flat terrain), it is not necessarily true in the case of complicated terrain. The alternative approach -- nonparametric methods -- estimate  distributions entirely based on data. Although robust and easy to implement, the accuracy of non-parametric methods strongly depend on the availability of data samples. In this paper, a nonparametric method, kernel density estimation (KDE), is used to describe the stochasticity of the turbulence. It is worth noting that the applicability of the framework developed in this context should not depend on the choice of density estimation techniques.}. To this end, an efficient method that leverages the huge amount of meteorological data that is readily available will be very useful for various aspects of wind turbine analysis. This paper describes such a technique. The contributions of this work are as follows: 1) we formulate a mathematical framework that is able to accurately represent wind flow using a low-complexity yet realistic stochastic model; 2) we implement a computational framework based on the mathematical framework that extracts statistical information from meteorology data and constructs an easy-to-use model; 3) we test the model by generating synthetic stochastic wind flow and comparing several wind turbulence statistics between synthetic and original flows.

In the following sections, section~\ref{math} focuses on constructing a low-complexity model of the random wind flow. Bi-orthogonal Decomposition (BD)~\cite{Venturi2008,Venturi2011} and Karhunen-Lo\`{e}ve expansion (KLE)~\cite{Ghanem1998} are used in constructing the low-dimensional wind model. A numerical example uses CWEX-11 (Crop Wind-energy EXperiment 2011) data is given in section~\ref{numerical example}. The algorithmic details of the computational framework is give in section~\ref{algorithm and implementation}. The result of the analysis is shown in section~\ref{results}. In order to quantitatively check the accuracy of the decomposition, several statistical properties of the synthetic stochastic wind flow and the original wind flow are compared. Finally, section~\ref{conclusion} concludes the paper with the main findings and their implications.

\section{Mathematical Framework}\label{math}
A stochastic wind can be denoted by $\bm v(\bm x,t,\bm\xi)$, where $\bm x=(x_1,x_2,x_3)=(x,y,z)$ is the three-dimensional coordinates with $x$, $y$, and $z$ represents along-wind, transverse, and vertical axis respectively, $t$ denotes time, and $\bm \xi$ represents stochastic variability. Wind velocity $\bm v=(v_1,v_2,v_3)=(u,v,w)$ consists of three wind speed components along the three spatial axes. 

The analysis is clearer when the mean component $\bar{\bm v}(\bm x)$ of the velocity data is removed so that only the fluctuation components $\bm u(\bm x,t,\bm\xi)$ are left, i.e.
\begin{equation}
 \bm u(\bm x,t,\bm\xi) = \bm v(\bm x,t,\bm\xi) - \bar{\bm v}(\bm x), \label{minus mean}.
\end{equation}
The mean component is defined as
\begin{equation}
 \bar{\bm v}(\bm x) = \frac{1}{|T|}\int_T \langle \bm v(\bm x,t,\bm\xi) \rangle dt. \label{mean}
\end{equation}
where $\langle \cdot \rangle$ denotes the average in the stochastic domain, $|T|$ is the span of the temporal domain, $\bar{\bm v}$ is the ensemble average.

The goal of this work is to construct a simple stochastic model that encodes temporal and spatial covariance, and preserve all the statistical information, such as spatial coherence and wind speed power spectral density (PSD), present in the collected meteorological data. We look for a model that has the form
\begin{equation}
\bm u =\sum K_{ijk} \bm X_i(\bm x)\ T_j(t)\ \bm\xi_k
\label{Eqn:First}
\end{equation}
where $K_{ijk}$ are coefficients, the deterministic functions $\bm X_i$ track spatial correlations, $T_j$ track temporal correlations and the random variables $\bm\xi_k$ encode the inherent variability. {\it The goal becomes to find the simplest possible representation that still encodes all the required information}. 

The main idea of the next section is to formulate a mathematical strategy of representing this meteorology data in terms of the smallest possible number of terms in Eqn.~\ref{Eqn:First}, by optimally designing $\bm X$, $T$, and $\bm\xi$. Noting that the data contains spatial, temporal and stochastic variabilities, we solve this problem in two stages. In the first stage, we decompose the data into temporal ($T(t)$) and coupled spatial-stochastic ($\bm \varPhi(\bm x,\bm\xi)$) parts through the concept of Bi-orthogonal Decomposition; in the second part we decompose the spatial-stochastic part into spatial ($\bm X(\bm x)$) and stochastic ($\bm\xi$) components using the concept of the Karhunen-Lo\`{e}ve decomposition. 

\subsection{Stage 1: A low-dimensional representation via Bio-orthogonal decomposition}\label{representation}

In the first stage, we are looking for a minimal representation of the data in the form
\begin{equation}
 \bm u(\bm x,t,\bm\xi) \approx \sum_{i=1}^{M} K_{i} \bm \varPhi_i(\bm x,\bm\xi) \ T_i(t).  \label{BD}
\end{equation}
where $\bm \varPhi_i(\bm x,\bm\xi)$ are stochastic spatial modes and $T_i(t)$ are temporal modes.\footnote{A general representation will be of the form
$\bm u(\bm x,t,\bm\xi)=\sum_{i,j} K_{ij} \bm \varPhi_i(\bm x,\bm\xi) \ T_j(t). $
where $i$ and $j$ are independent indices. Consequently, the expression on the right hand side has an dramatically large number of terms. This is handled utilizing the Schmidt decomposition theorem~\cite{Schmidt1907}, which states that -- any representation of a tensor product space $\mathscr{H} = \mathscr{H}_1 \otimes \mathscr{H}_2$ can be expressed as linear combination of tensor product of basis functions $\bm \varPhi_i \otimes \varPsi_i$, where $\bm \varPhi_i \in \mathscr{H}_1$, $\varPsi_i \in \mathscr{H}_2$. As a result, the representation can be reduced to Eqn.~\ref{BD}.}
Our goal becomes searching for the best choices for $\bm \varPhi_i$ and $T_i$ such that the decomposition uses the least number of terms M that will give us an accurate representation of $\bm u(\bm x,t,\bm\xi)$. We pose this as an optimization problem. To do so, we define the error in this representation and design $\bm \varPhi_i$ and $T_i$ that minimize this error. The error, denoted by $\varepsilon$, is defined as the low-complexity stochastic model subtracted from the true data 
\begin{equation}
 \bm \varepsilon(\bm x,t,\bm\xi) = \bm u(\bm x,t,\bm\xi) - \sum_{i=1}^M K_i \bm \varPhi_i(\bm x,\bm\xi) \ T_i(t).  \label{error}
\end{equation}
Approximation theory suggest that the best choice for the functions $\bm \varPhi_i$ and $T_i$ (we will also call them modes) is when they are orthogonal to each other~\cite{Hunter2005}. Thus, we set $T_i,~i=1,\ldots,M$ to be orthogonal to each other in the time domain and $\bm \varPhi_i, ~i=1,\ldots,M$ to be weakly orthogonal in spatial domain. Mathematically, this is denoted in terms of the inner products:
\begin{equation}
 \langle T_i,T_j \rangle_T = \int_T T_i(t) T_j(t)dt = \delta_{ij} \label{temporal inner}
\end{equation}
and 
\begin{equation}
 \langle \bm \varPhi_i,\bm \varPhi_j \rangle_X = \int_{\bm X} \overline{\bm \varPhi_i} \cdot \overline{\bm \varPhi_j} \ d\bm x = \delta_{ij}, \label{spatial inner}
\end{equation}
where $\overline{\bm \varPhi_i}$ denotes the expectation of the spatial-stochastic mode, i.e.
\begin{equation}
 \overline{\bm \varPhi_i(\bm x)} = \int \bm \varPhi_i(\bm x,\bm\xi) W(\bm\xi) \ d\bm\xi, \label{spatial inner 2}
\end{equation}
and $W(\bm\xi)$ is the multivariate joint probability density of random variables in the set $\bm\xi$. 

Note that the error $\bm \varepsilon$ itself is an random field. We construct an associated scalar value with this random field to accomplish subsequent optimization. To this end, an error functional $\mathscr{E}$ is defined as the norm of $\bm \varepsilon$, i.e.
\begin{equation}
 \mathscr{E} = \int_T \langle \bm \varepsilon, \bm \varepsilon \rangle_X dt.  \label{error functional}
\end{equation}
The error-functional is simply the inner product (i.e. an average) of the field over space, time and stochastic dimensions. Note that the error functional depends on the choice of functions $T_i(t)$ and $\bm \varPhi_i(\bm x,\bm\xi)$, i.e. $\mathscr{E}[T_1,\cdots,T_M,\bm \varPhi_1,\cdots,\bm \varPhi_M]$. 

We now search for temporal functions that minimizes this error functional. This is accomplished by applying the calculus of variations and solving the associated Euler-Lagrange equations. We provide full details of the derivation in Appendix~\ref{reperror}. This optimization reduces to the solution of an eigenvalue problem, where the eigenfunctions give the temporal modes:
\begin{equation}
 \mu_i T_i(t) = \int_T C(t,t') T_i(t')dt',  \label{BD eigen problem}
\end{equation}
where $C(t,t')$ is the temporal covariance matrix of the meteorology data, and $\mu_i$ are the eigenvalues of this temporal covariance matrix. $C(t,t')$ is constructed by taking the inner product of the data in the spatial domain, i.e
\begin{equation}
 C(t,t') = \langle \bm u(\bm x,t,\bm\xi),\bm u(\bm x,t',\bm\xi) \rangle_X.  \label{BD autocorr}
\end{equation}
Once $\mu_i$ and $T_i$ are solved for, the spatial-stochastic functions are calculated using 
\begin{equation}
 \bm \varPhi_i(\bm x,\bm\xi)=\frac{1}{\sqrt{\mu}_i}\langle \bm u(\bm x,t,\bm\xi),T_i(t) \rangle_T. \label{best Phi}
\end{equation}
The first $M$ (usually $M \sim 3-6$) eigenvalues and eigenfunctions usually represent the data exceedingly well~\cite{Ganapathysubramanian2007}. Thus, the first stage of the decomposition of the wind data results in representation involving temporal functions $T_i$ and spatial-stochastic functions $\bm \varPhi(\bm x,\bm\xi)$. Eqn.~\ref{BD} now becomes
\begin{equation}
 \bm u(\bm x,t,\bm\xi) = \sum_{i=1}^M \bm a_i(\bm x,\bm\xi) \ T_i(t),  \label{new BD}
\end{equation}
where 
\begin{equation}
 \bm a_i(\bm x,\bm\xi) = K_i\bm \varPhi_i(\bm x,\bm\xi) = \sqrt{\mu_i} \ \bm \varPhi_i(\bm x,\bm\xi).   \label{spatial stochastic modes}
\end{equation}
are spatial-stochastic modes. 
The second stage of the framework involves decomposing the spatial-stochastic functions $\bm a_i(\bm x,\bm\xi)$ into spatial functions $\bm X(\bm x)$ and independent random variables $\bm\xi$.

\newtheorem{remark}{Remark}
\begin{remark}
We chose to decompose the data into spatial-stochastic and temporal parts in the first stage of the decomposition. This decomposition is one of three possible decomposition choices. These choices are enumerated in Table~\ref{choices} where $\mathbb X$, $\mathbb T$, and $\varSigma$ denote the spatial, temporal, and stochastic domains respectively. Such decompositions have been explored in other works. For instance, Venturi et al.~\cite{Venturi2008,Venturi2011} investigated Type 1 decomposition. The decomposition suggested by Type 2 can be achieved by using generalized polynomial chaos~\cite{Xiu2002,Xiu2003}. Mathelin et al.~\cite{Mathelin2009} modeled uncertain cylinder wake using a Type 3 decomposition. It can be shown that Type 1 and Type 3 decompositions result in identical results. \hfill \qed
\begin{table}[htbp]%
\caption{Choices for Bi-orthogonal Decomposition.}
\centering
\begin{tabular}{ccc}
\hline %
Type&Space 1&Space 2\\\hline%
1&$\mathbb X \times \varSigma$&$\mathbb T$\\
2&$\mathbb X \times \mathbb T$&$\varSigma$\\
3&$\mathbb T \times \varSigma$&$\mathbb X$\\\hline 
\end{tabular}
\label{choices}
\end{table}
\end{remark}

\begin{remark}
The choice of the inner products (temporal Eqn.~\ref{temporal inner} and spatial-stochastic Eqn.~\ref{spatial inner}) affect the properties of the decomposition. In particular, there are several different ways in which we can define an inner product over the spatial-stochastic modes (i.e. different ways to average over space and stochastic dimensions). These include the following possibilities: 
\begin{equation}
\langle \bm \varPhi_i,\bm \varPhi_j \rangle_0 = \int_{\bm X} \overline{\bm \varPhi_i} \cdot \overline{\bm \varPhi_j} \ d\bm x,\nonumber
\end{equation}
\begin{equation}
\langle \bm \varPhi_i,\bm \varPhi_j \rangle_1 = \int_{\bm X} \overline{\bm \varPhi_i \cdot \bm \varPhi_j} \ d\bm x,\nonumber
\end{equation}
\begin{equation}
\langle \bm \varPhi_i,\bm \varPhi_j \rangle_2 = \int_{\bm X} \overline{\bm \varPhi_i \cdot \bm \varPhi_j} - \overline{\bm \varPhi_i} \cdot \overline{\bm \varPhi_j}
  \ d\bm x,\nonumber
\end{equation}
The first inner-product (denoted by $\langle \cdot \rangle_0$) is a spatial integral of the product of expected values, while the second inner-product (denoted by $\langle \cdot \rangle_1$ ) is the expectation of the spatial integral.
It can be shown that by taking inner products, $\langle \cdot \rangle_h$, of type $h=0,1,2$, we obtain optimal representations with respect to mean, second-order moment, and standard deviation of the data, respectively. We have chosen to focus on representation that is optimal in the mean sense~\cite{Venturi2008}. \hfill \qed
\label{remark inner}
\end{remark}

\subsection{Stage 2: Karhunen-Lo\`{e}ve expansion of spatial stochastic modes}\label{KLE section}
In this stage, we decompose the spatial-stochastic functions $\bm a_i(\bm x, \bm\xi)$ into a spatial part and a set of uncorrelated random variables. First, the mean of spatial-stochastic functions $\bar{\bm a}(\bm x)$ are removed so that only the fluctuation component is left for analysis.
\begin{equation}
 \bm \alpha(\bm x,\bm\xi) = \bm a(\bm x, \bm\xi) -  \bar{\bm a}(\bm x). \label{isolate mean KLE}
\end{equation}
Following the rational of Bi-orthogonal Decomposition (Eqn.~\ref{BD}), our goal is to decompose $\bm \alpha(\bm x,\bm\xi) $ into a minimal set of linear combination of deterministic spatial functions and uncorrelated random variables, that is
\begin{equation}
 \bm \alpha(\bm x,\bm\xi) \approx \sum^N_{i=1} \, C_i\, \xi_i\, \bm X_i(\bm x).  \label{approx KLE}
\end{equation}
where $C_i$ are coefficients of the expansion, $\xi_i$ is a set of uncorrelated random variables, $\bm X_i$ are deterministic spatial functions. We pose this decomposition problem as an optimization problem, where the optimization problem is to minimize the error. The representation error is defined as
\begin{equation}
 \bm \varepsilon = \bm \alpha(\bm x, \bm\xi) - \sum^N_{i=1} \, C_i\, \xi_i\, \bm X_i(\bm x)  \label{rep error KLE}
\end{equation}
This is a standard formulation of the Karhunen-Lo\`{e}ve expansion. The goal of KLE is to find the optimal choice for functions $\bm X_i$ such that the representation error is minimized with a finite number ($N$) of expansion terms. We briefly describe the mathematical framework of KLE below~\cite{Ghanem1998}: The representation error is converted into a cost-functional for optimization by simply considering the mean-square error (i.e. the inner-product)
\begin{equation}
 \mathscr{E}^2 = \int_{\bm X}\, \overline{\bm \varepsilon}^2\,d\bm x.  \label{error functional KLE}
\end{equation}
Minimization of the error-functional results in an eigenvalue problem, whose eigenfunctions are the desired spatial functions $\bm X_i$
\begin{equation}
  \int_{\bm X} R(\bm x_1,\bm x_2) \, \bm X_{i}(\bm x_2) \, d\bm x_1 = C^2_i \, \bm X_{i}(\bm x_2),
  \label{EP_KLE}
\end{equation}
where $R(\bm x_1,\bm x_2)$ is the covariance kernel constructed from the spatial-stochastic functions $\bm \alpha(\bm x, \bm\xi)$. We denote $C^2_i = \lambda_i$ so that $\lambda_i$  and $\bm X_i(\bm x)$ are the eigenvalues and  the eigenvectors of the covariance kernel. Eqn.~\ref{approx KLE} becomes
\begin{equation}
 \bm \alpha(\bm x,\bm\xi) \approx \sum^N_{i=1} \, \sqrt{\lambda_i}\, \xi_i\, \bm X_i(\bm x).  \label{KLE}
\end{equation}
The final stage of the formulation is to identify the probability distributions of the uncorrelated random variables $\xi_i$. Note that we have the same number of realizations of $\xi_i$ as the number of random samples of the stochastic wind. Each realization is computed by inverting Eqn.~\ref{KLE} for $\xi_i$:
\begin{equation}
 \xi_i = \frac{1}{\sqrt{\lambda_i}} \int_{\bm X} \bm \alpha \bm X_i(\bm x) d\bm x.  \label{XI}
\end{equation}
The choice of technique to construct the probability distribution of $\xi_i$ given a finite number of observations of $\xi_i$ is crucial. We utilize Kernel Density Estimation (KDE) methods~\cite{Scott1992} to construct the probability distributions of $\xi_i$ in a non-parametric way.

\begin{remark}
In order to solve the eigen-problem (Eqn.~\ref{EP_KLE}), the covariance kernel must first be calculated from the data. The Wiener-Khinchin theorem allows computing the covariance in a very efficient way. Representing the spatial-stochastic data as a matrix $I$, the covariance, $R$, is computed as: 
\begin{equation}
  R = \mathscr{F}^{-1}(\mathscr{F}(I) \times \mathscr{F}(I)').
\end{equation}
where $\mathscr{F}(I)'$ is the complex conjugate of Fourier transform of I, and the diagonal entries of $R$ contains the covariance. Numerical details for solving this generalized eigenvalue problem is included in Appendix~\ref{numericaleigen} and a software available~\cite{CODELINK}. \hfill \qed
\end{remark}

\begin{remark}
In KDE, the PDF of a random variable $\xi$ is estimated as
\begin{equation}
 p(\xi)=\frac{1}{Nh}\sum_{i=1}^N K \left( \frac{ \xi-\xi_i }{h} \right),  \label{KDE}
\end{equation}
where $K(\cdot)$ is the kernel which is a symmetric function that integrates to one, and $h>0$ is a smoothing parameter called the bandwidth. Common choices for $K(\cdot)$ are the multivariate Gaussian density function. The choice of the bandwidth is critical since a small values of $h$ result in the estimated density with many 'wiggles', while large values of $h$ result in very smooth estimations that do not represent the local distributions. Several approaches have been developed to chose proper $h$ values. A good review on bandwidth selection can be found in~\cite{Turlach1993}. We advocate using a simple formula based on Silverman's rule~\cite{Silverman1986}. The optimal choice for $h$ according to Silverman's rule is given by
\begin{equation}
 h=\left( \frac{4\hat{\sigma}}{3N} \right) ^ {\frac{1}{5}} \approx 1.06 \hat{\sigma} n^{-\frac{1}{5}},  \label{silverman}
\end{equation}
where $\hat{\sigma}$ is the standard deviation of the samples. \hfill \qed
\end{remark}

Based on the above three techniques, we are now able to construct a space-time decomposition of the wind field snapshots and preserve the spatial and temporal correlations. In next section, we look at a numerical example that illustrates the power of this framework.

\section{Numerical Example: CWEX-11}\label{numerical example}
We focus on constructing realistic low-complexity models from experimental data. We illustrate the methodology based on a recent meteorological experiment named CEWX-11. In this context, the CWEX-11 data is ideal to test the framework and demonstrate the advantages of the framework. It is straightforward to extend this analysis to any other meteorological dataset that expands the scale of measurements of the CWEX-11 experiment. \footnote{The CWEX-11 data measures the wind profile at two locations: 4.5 meters and 10 meters. Obviously, this data does not describe the wind speed profile at hub height. Additionally, measurements at only two different heights are certainly not enough to capture all the spatial stochasticity in the wind.}

CWEX-11 is a collaborative experiment (Iowa State University (ISU) and the University of Colorado (UC), assisted by the National Center for Atmospheric Research (NCAR)). CWEX-11 and its 2010 counterpart (CWEX-10) in the Crop Wind-energy EXperiment, address observational evidence for the interaction between large wind farms, crop agriculture, and surface-layer, boundary-layer, and mesoscale meteorology (Rajewski et al. 2012)~\cite{CWEX-BAMS}. In the experiment, a surface flux station was installed to the south of a wind turbine, which makes the station measure upstream inflow of the wind turbine due to the fact that predominant summer winds in Iowa originate from south to slightly south-east. The surface station was equipped with a CSAT3 sonic anemometer that was located at height 4.5m and an RMY propeller and vane anemometer that was located at height 10.0m. The former measured wind speed in 3 directions at 20 Hz whereas the latter gave wind speed amplitude and its direction at 1 Hz. The experiment started 
on June 29 and lasted for 48 days.\footnote{Excluding the days that mostly have opposite wind direction and the days with sudden rain, only data for 28 days that had the desired weather condition and wind direction were used in the analysis.}

\begin{figure}[htbp]
\centering
\includegraphics[width=0.8\textwidth]{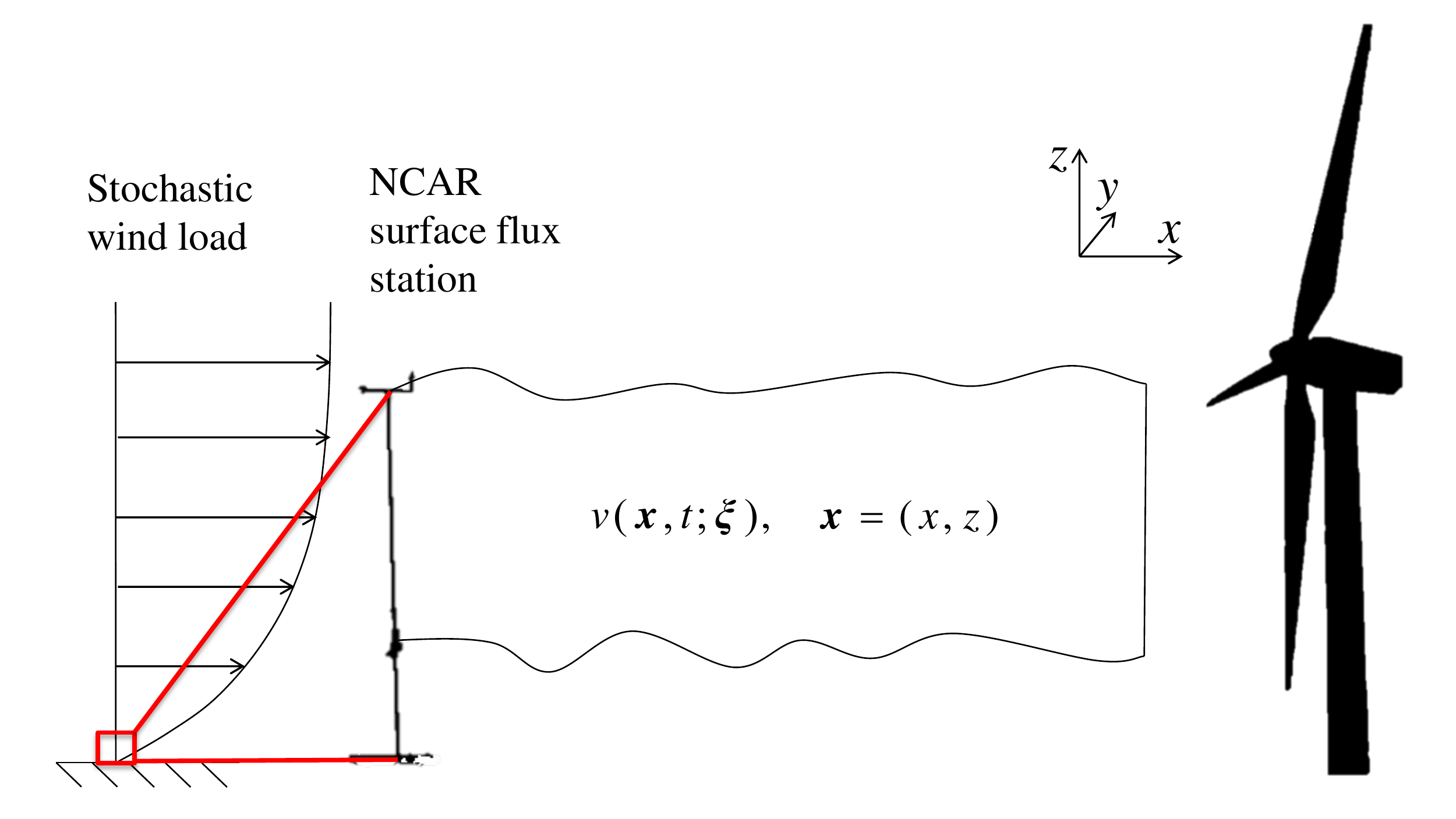}
\caption{Setup of experiment}
\label{expsetup}
\end{figure}

A schematic of the experiment is shown in Fig.~\ref{expsetup}. Note that the figure is for illustrative purpose and is not drawn to scale. In this experiment, only wind speed magnitude was analyzed, but it is straightforward to perform analysis on any one of the three components of turbulence. 

Ideally, we would like to have multiple wind field snapshots taken at regular intervals on a vertical plane that is perpendicular to the rotor. However, due to the constrains of the experiment, only two time-series were measured at heights 4.5m and 10.0m. Linear interpolation of the two time-series at 18 spatial points on $z$ direction were performed to simulate more vertical measurements. The wind snapshots are constructed by using {\it Taylor's frozen turbulence hypothesis}~\cite{taylor1938}. Specifically, all measurements in the resulting 20 time-series over certain interval were treated as taken at the same instant in the interval. This  provides us along-wind ($x$ direction) measurements. Finally, one 2D snapshot taken at certain instant is constructed\footnote{Since no measurements were taken in the transverse direction, only 2D wind fields are analyzed in this work, though the mathematics is agnostic to dimensionality}. Similarly, snapshots of the wind field at other instants were created to represent one full day of measurements. We have data for 28 
such days. The full meteorology data curated into this form can be represented as following matrix
\[\begin{array}{cccc}
\bm v_{11}(x,z) & \bm v_{12}(x,z) & \cdots & \bm v_{1m}(x,z)\\
\bm v_{21}(x,z) & \bm v_{22}(x,z) & \cdots & \bm v_{2m}(x,z)\\
\vdots & \vdots & \  & \vdots\\
\bm v_{n1}(x,z) & \bm v_{n2}(x,z) & \cdots & \bm v_{nm}(x,z) \end{array}\]
In the matrix, each element $\{\bm v_{ij}(x,z),\ i=1,2,\cdots,m; \ j=1,2,\cdots,n\}$ represents a snapshot of wind field. Each column contains snapshots that span an entire day. Each row represents data for a particular interval measured over different days. Thus, each row of data can be considered to be {\it realizations} of the stochastic wind field at a particular interval. 

It is noteworthy that the length of time interval that was used to obtain wind field snapshots can be arbitrarily chosen. Although, most guidelines for wind turbine design and siting suggest considering wind variabilities over a period of 10 minutes for wind data analysis~\cite{GermanischerLloyd2010}, in this paper, analysis based on different choices of time intervals were performed and compared. In the results section, we will show that using 10 minutes time interval to construct wind field snapshots is an reasonable choice. 

\section{Algorithm \& Implementation}\label{algorithm and implementation}
The algorithmic details of the framework is outlined in Table \ref{steps}. The meteorological data is provided in the common NetCDF (Network Common Data Form) format. The original meteorological data contains a variety of information, including wind speed and orientation, temperature, humidity, surface $\mathrm{CO}_2$ flux. For this analysis we only consider the wind speed data. The wind speed data is first extracted from the NetCDF file using MATLAB. The matrix of wind flow snapshots is then constructed. The temporal covariance matrix is calculated from the snapshot matrix and stored as a datafile. Next, the eigenvalue problem (corresponding to the Bi-orthogonal Decomposition) is solved to get temporal and spatial-stochastic modes. Following this, the spatial-stochastic modes are decomposed into spatial functions and uncorrelated random variables via the Karhunen-Lo\`{e}ve expansion. This is accomplished by solving an eigenvalue  problem. The eigenvectors computed are the spatial functions and they are used 
to compute realizations of the uncorrelated random variables (via Eqn.~\ref{XI}). Finally, a Kernel Density Estimator is used to construct the PDFs of random variables $\xi^i_j$. At this time, all the necessary elements for constructing synthetic wind flow were known. Finally, synthetic wind realizations are constructed by reversing above process. The complete framework is implemented in C++. SLEPc \cite{Hernandez2005}, which is based on PETSc \cite{petsc-web-page}, was used in solving eigenvalue problems. 

\begin{table}[htbp]
\caption{Steps of computational framework.}
\centering
\begin{tabular}{p{0.8\textwidth}}
\hline%
\begin{enumerate}
 \item Data preparation
  \begin{enumerate}
   \item Extract wind speed data from NetCDF files
   \item Construct matrix of wind speed field snapshots
  \end{enumerate}
 \item Bi-orthogonal Decomposition
  \begin{enumerate}
   \item Calculate temporal covariance from snapshots matrix
   \item Solve eigenvalue problem, get eigenvalues and eigenfunctions (temporal modes)
   \item Solve spatial stochastic modes
  \end{enumerate}
 \item Karhunen-Lo\`{e}ve expansion
  \begin{enumerate}
   \item Calculate covariance kernel of each spatial mode (M spatial modes in total)
   \item Solve eigenvalue problem for each spatial mode, get basis spatial functions
   \item Calculate observations of each random variable $\xi^i_j$ based on Eqn.~\ref{XI}
  \end{enumerate}
 \item Kernel density estimation
  \begin{enumerate}
   \item Estimate PDF for each random variable $\xi^i_j$ based on observations
  \end{enumerate}
 \item Construct synthetic wind flow
  \begin{enumerate}
   \item Get random samples from random variables $\xi^i_j$
   \item Construct stochastic spatial modes according to Eqn.~\ref{KLE}
   \item Construct synthetic wind flow using Eqn.~\ref{BD}
  \end{enumerate}
\end{enumerate} \\\hline
\end{tabular}
\label{steps}
\end{table}

\section{Results}\label{results}
As stated earlier, we use the meteorological data from the CWEX-11 experiment. The data consists of wind speed amplitude and direction measure at 1 Hz over 28 days. We will describe the stage wise construction of the low-dimensional model and its accuracy in the next few subsections. 

\subsection{BD results}\label{BD results}

\begin{figure}[htbp]
 \centering
 \includegraphics[width=0.8\textwidth]{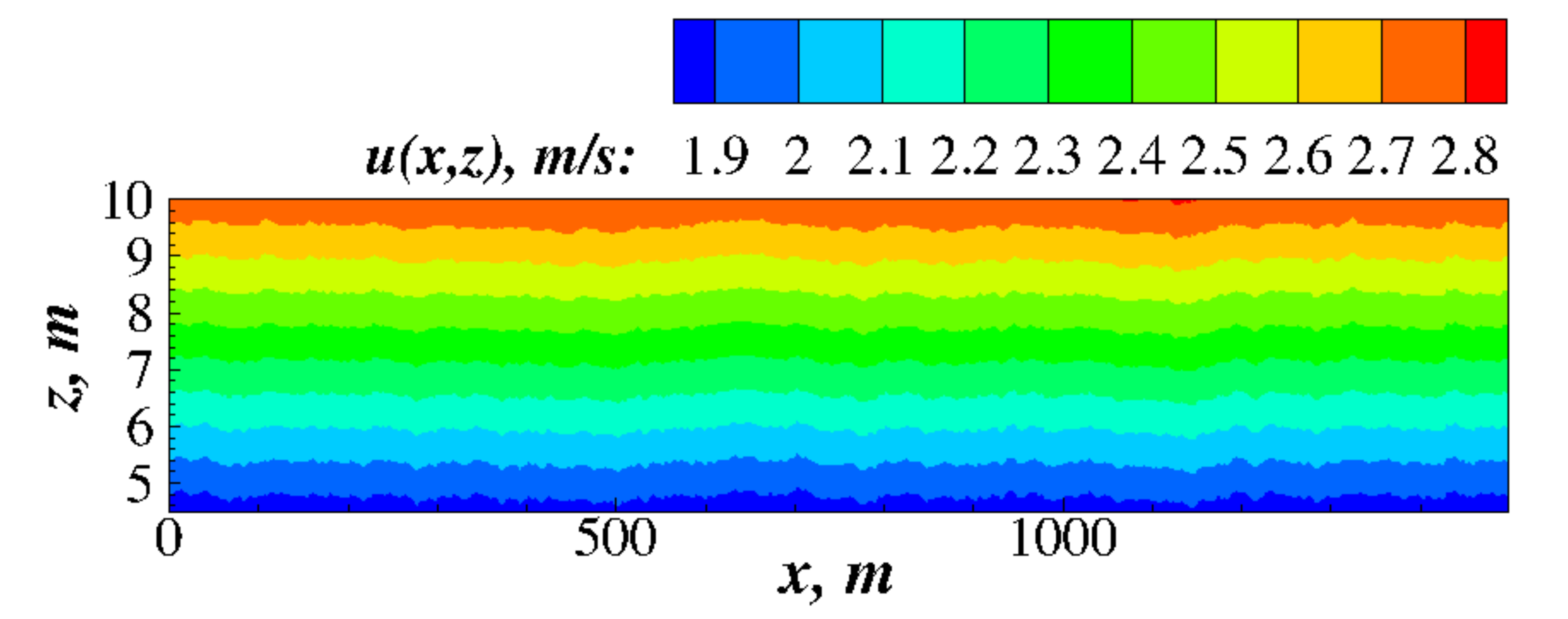}
\caption{10 minute period, averaged across all 144 time intervals and 28 realizations}\label{meanplot}
\end{figure}

As discussed previously, snapshots curated based on 10 minutes interval are used first. As a result, there are 144 snapshots in one day and we have 28 such days (samples). The mean component of the flow that is averaged across all 144 snapshots and 28 realizations according to Eqn.~\ref{mean} is shown in Fig.~\ref{meanplot}. This plot shows us the mean wind profile, and clearly illustrates the vertical wind shear. Note that the x-axis represents length of the snapshots that are calculated based on 10 minute interval and average wind speed.

\begin{figure}[htbp]
\centering
\includegraphics[width=0.7\textwidth]{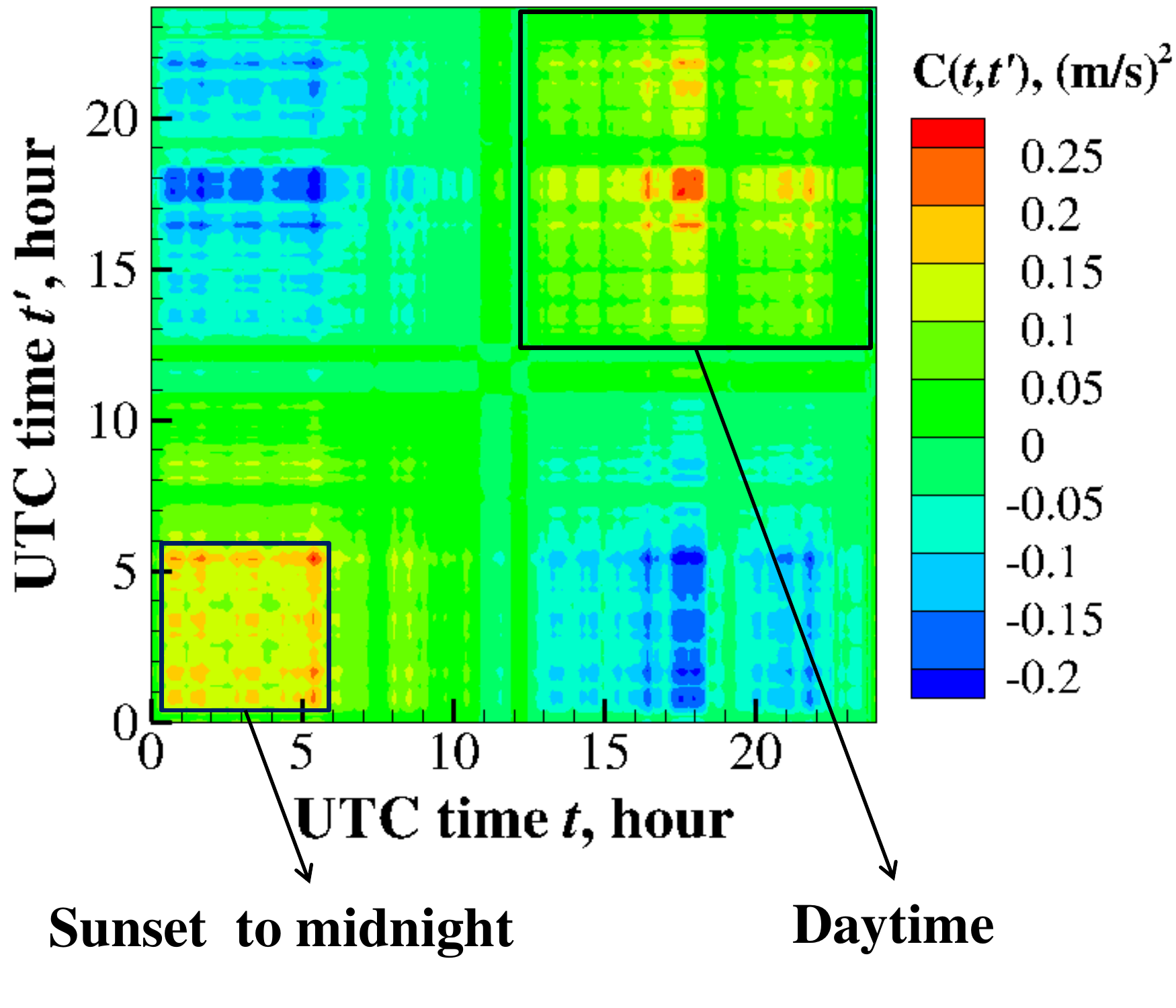}
\caption{Covariance function $C(t,t')$}
\label{tempcor}
\end{figure}

\begin{figure}[htbp]
\centering
\includegraphics[width=0.95\textwidth]{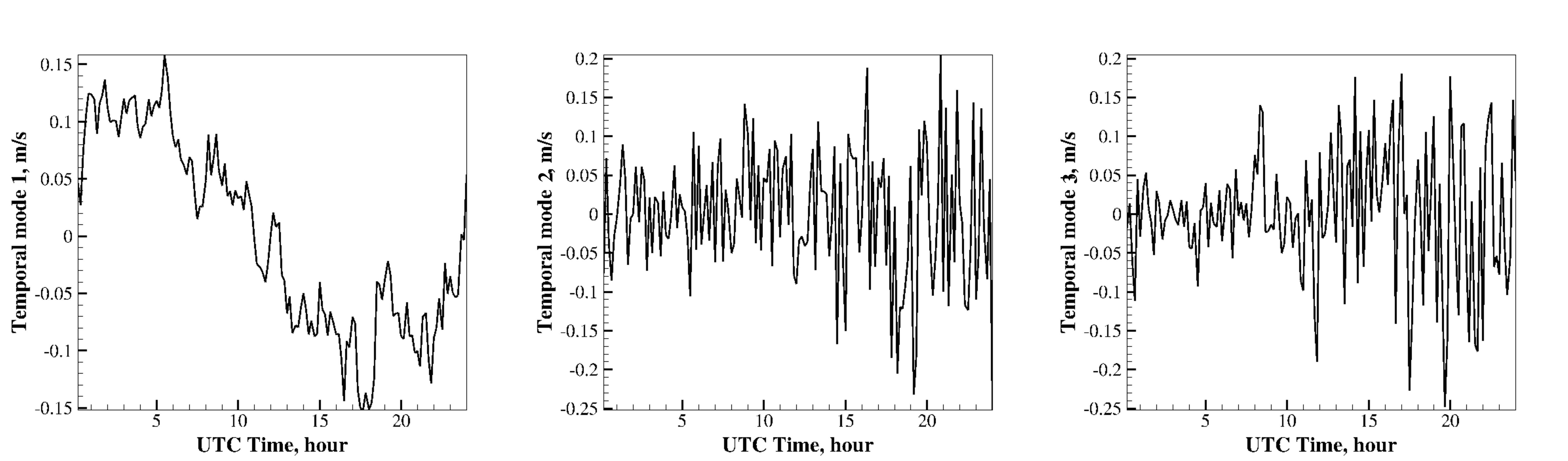}
\caption{First three eigenvectors of temporal covariance function \label{eigenfunc}}
\end{figure}

The temporal covariance are calculated based on Eqn.~\ref{BD autocorr}. By taking the inner products across the spatial domain, we construct the covariance function $C(t,t')$ which is shown in Fig.~\ref{tempcor}. Note that $C$ has block structure, with regions of high covariance (marked by the solid boxes) along the diagonal and regions of large negative covariance along the off-diagonal. This structure of the covariance function follows the dynamics of stable, and unstable stratification of the atmospheric boundary layer seen in the US central plains. We discuss this by dividing the analysis into several distinct time periods (marked by the solid boxes). Following meteorological practice, the data starts at Coordinated Universal Time (UTC, or Greenwich time) 00:00. The first period is UTC 00:00-06:00 (CST 18:00-00:00) which corresponds to the time between sunset to midnight. In this period, insolation and, thus, heating is gradually cut off and the temperature of atmosphere cools down (from the ground up) 
due to rapid cooling of the ground. Because of the heavier density of the cooler air, the cooler air stays at the bottom (close to the surface). This generates a stably stratified boundary layer that does not change during the duration of the night. This results in the high covariance between adjacent time periods marked in the lower left box. The second period is the region of reduced covariance between UTC 06:00-13:00 which is basically the time from midnight to shortly after sunrise. During summer, the nocturnal Great Plains low-level jet (LLJ) exists in Iowa. For a good overview of LLJ, see Jiang et al.~\cite{Jiang2007}. The LLJ causes low level turbulence and enhanced mixing which reduces the covariance between adjacent-time wind fields. The third time period is between UTC 13:00-00:00 which corresponds to sunrise and day time. Because of the sunrise, the ground is rapidly heated. As a result, the warmer air near the ground becomes buoyant and rises rapidly with its place being taken by compensating 
flow of colder air from higher up in the atmospheric boundary layer. Eddies are thereby created, changing from smaller scale to larger scales, finally becoming the circulatory motion that crosses the entire boundary layer so that the high speed free stream flow is brought in the circulation. This phenomenon usually happens at noon, which can also be found in Fig.~\ref{eigenfunc} at UTC 18:00. In addition, the covariance function exhibits periodic behavior that is caused by the periodical motion of eddies. 

\begin{figure}[htbp]
\centering
\includegraphics[width=0.9\textwidth]{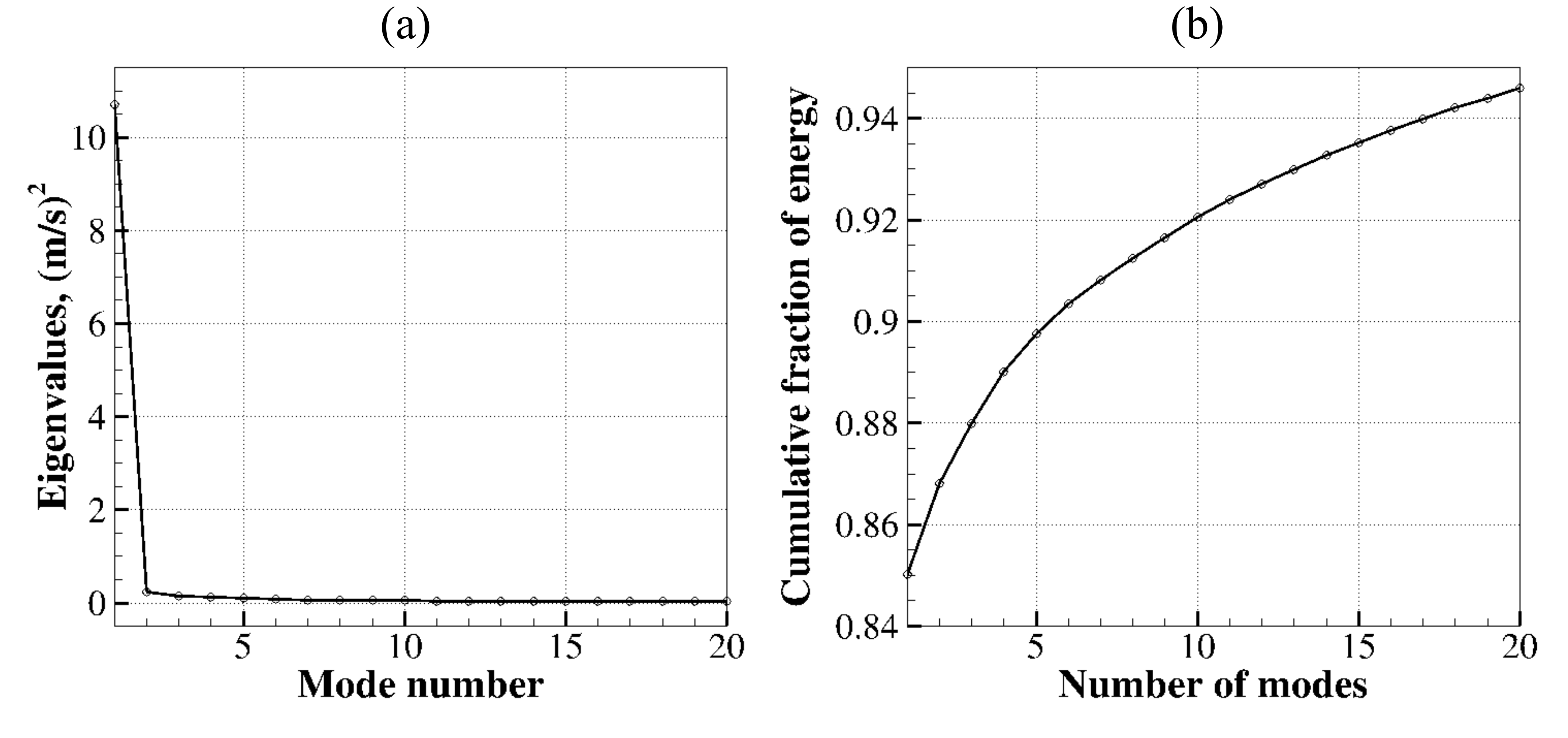}
\caption{Spectrum of $C(t,t')$\label{eigenvalue&sm}}
\end{figure}

\begin{figure}[htbp]
\centering
\includegraphics[width=0.6\textwidth]{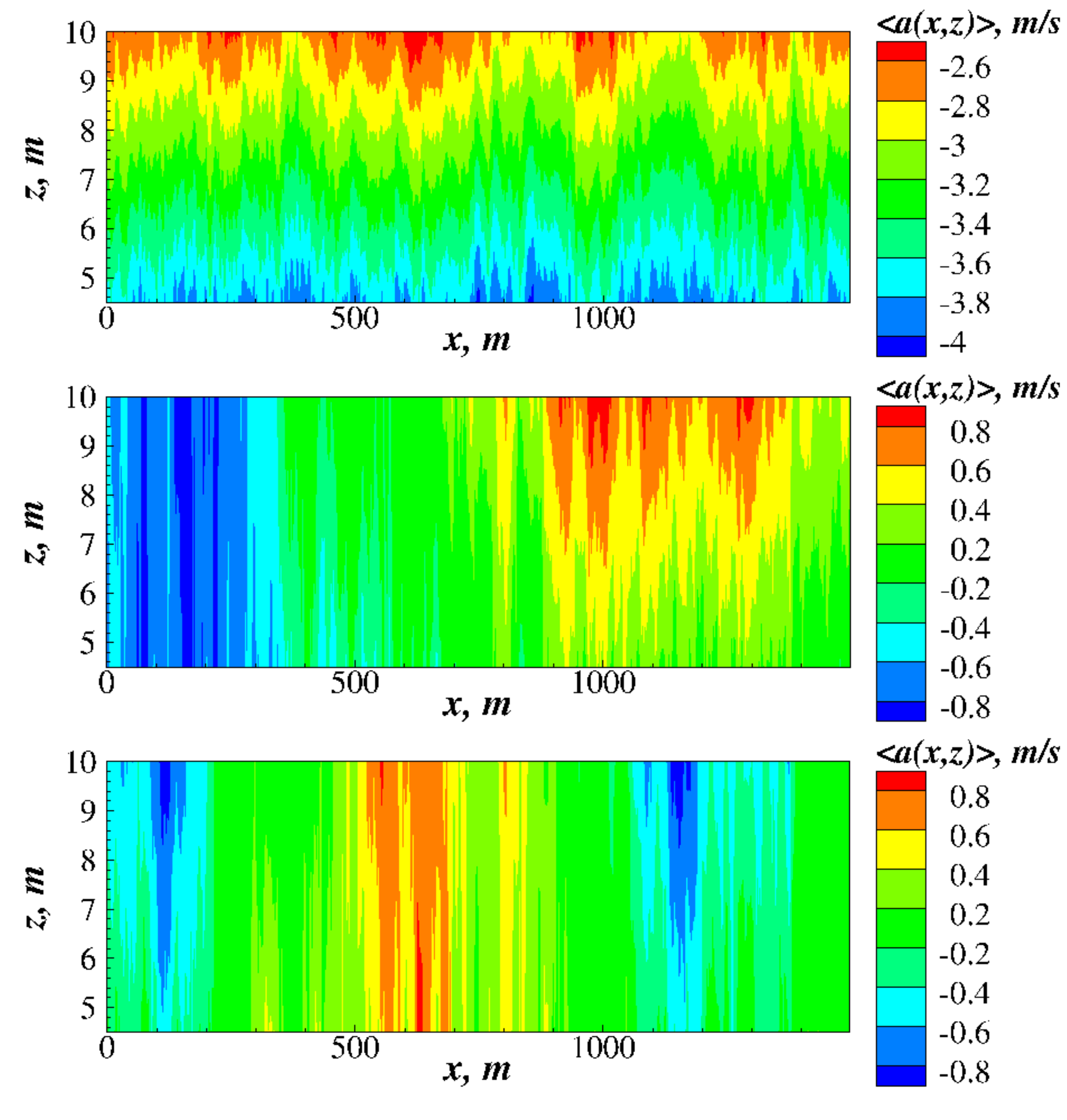}
\caption{Stochastic spatial modes of $C(t,t')$\label{spamodes}}
\end{figure}

The temporal covariance function plays a very important role in this analysis since it provides almost all the needed information that describes the behavior of the wind flow. Compared to the original meteorological data, the temporal covariance is much easier to store and transmit, yet contains nearly the same amount of information as the original meteorological data. The eigenvalues of the temporal covariance function are shown in Fig~\ref{eigenvalue&sm}. Fig~\ref{eigenvalue&sm}(a) compares the relative magnitudes of eigenvalues; the magnitudes of the eigenvalues provide a notion of how much energy about the data is stored in each spatial mode~\cite{LUMLEY}. Notice that the first eigenvalue is much larger than the other subsequent eigenvalues, which means that the first mode contains the largest portion of the energy in the turbulence field. In Fig~\ref{eigenvalue&sm}(b) this is represented as the cumulative fraction of energy contained in the first $k$ modes. This plot provides us with a precise notion of how many 
terms are needed to incorporate~\cite{LUMLEY}, say, 90\% of the information available in the data into the low-complexity model. The first five eigenvalues cover about 90 \% of the total energy of the turbulence field, which ensures that a five term decomposition will have a 90\% accuracy of representation. This illustrates the advantage of Bi-orthogonal Decomposition. As a reduced order model, Bi-orthogonal Decomposition is able to represent a random flow with much fewer random modes compared to the number of original snapshots. In other words, by using BD, we reduced the terms in representing the given stochastic random flow from 144 snapshots to a simple 5-variable parametrization. 
Fig.~\ref{eigenfunc} shows the first three eigenfunctions of $C(t,t')$. Note that the eigenfunctions are the temporal functions $T_i$. They track how the fluctuation component of the random flow varies during a day. Based on the plot, the largest mean (and standard deviation, which is not reported here) of wind speed occurs at UTC 18:00 because of the reason stated previously. From Fig.~\ref{eigenvalue&sm}, we notice that the first mode carries about 85\% of the total energy in the turbulence. Therefore, the trend of diurnal fluctuation of the turbulence field can be mostly seen in the first temporal mode, with the next two temporal modes look like white noise. 
Once eigenvalues and eigenvectors for temporal covariance function are solved, the spatial-stochastic modes $a_i(\bm x,\bm\xi)$ can be constructed (Eqn.~\ref{best Phi} and Eqn.~\ref{spatial stochastic modes}). Fig.~\ref{spamodes} shows the expectation of the first three spatial modes. It is clear that the first spatial mode that carries the largest part of turbulence energy describes the vertical shear pattern, while the second mode describes the lateral shear pattern of the wind field. Higher modes that account for more complicated turbulence structures are insignificant since they contain limited turbulence energy.

\subsection{KLE results}\label{KLE results}
The spatial-stochastic modes are decomposed into deterministic modes and random variables using KLE. KLE results of the first three spatial modes are shown in Table~\ref{tab:KLE}. Covariance matrices, eigenvalues, deterministic spatial modes, and probably distributions of random variables are shown in the table from top to bottom. As discussed above, the distributions of $\xi^i_j$ are estimated by using KDE. 

From Fig.~\ref{spamodes} we notice that as the mode order increases, the complexity of spatial mode also grows. As a result, higher order spatial modes require more KLE terms to be accurately represented. This statement can be verified by comparing results of the first three spatial modes. For instance, the smoothness of the covariance kernels, the complexity of deterministic modes, and the required numbers of KLE terms to achieve 90\% of representation accuracies all increases as the mode number becomes higher.

It is worth to mention that the PDFs of $\xi^i_j$ are not always Gaussian. Although Gaussian assumption is commonly used in turbulence analysis, non-Gaussian variables do exist. Notice that each random sample corresponds to one day wind history, the randomness of this long term stochastic process may be different from the randomness of the more commonly used 10-minute average wind speed whose distribution is always regarded as Gaussian or Weibull. 
\begin{sidewaystable}
  \begin{tabular}
   {ccc} \hline $a_1(\bm x,\bm \xi)$ & $a_2(\bm x,\bm \xi)$ & $a_3(\bm x,\bm \xi)$ \\
    \hline 
    \includegraphics[width=0.3\textwidth]{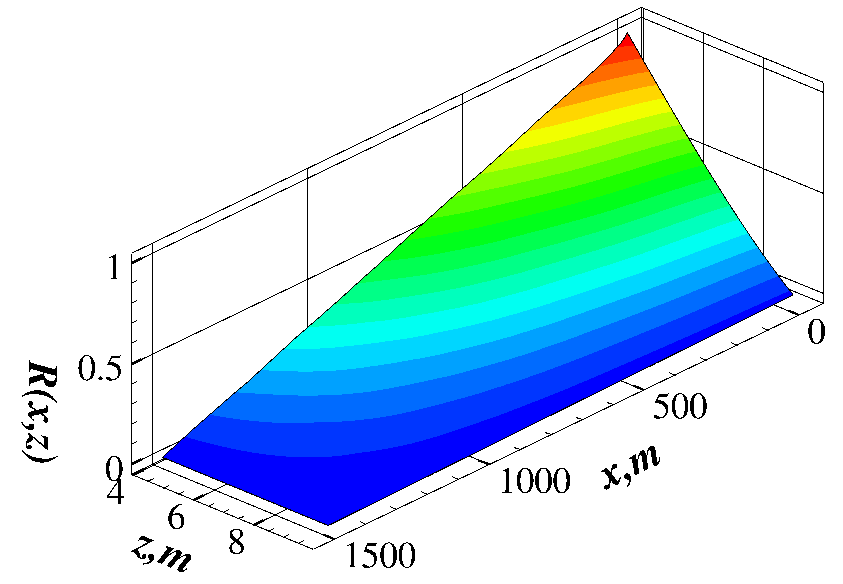} & 
    \includegraphics[width=0.3\textwidth]{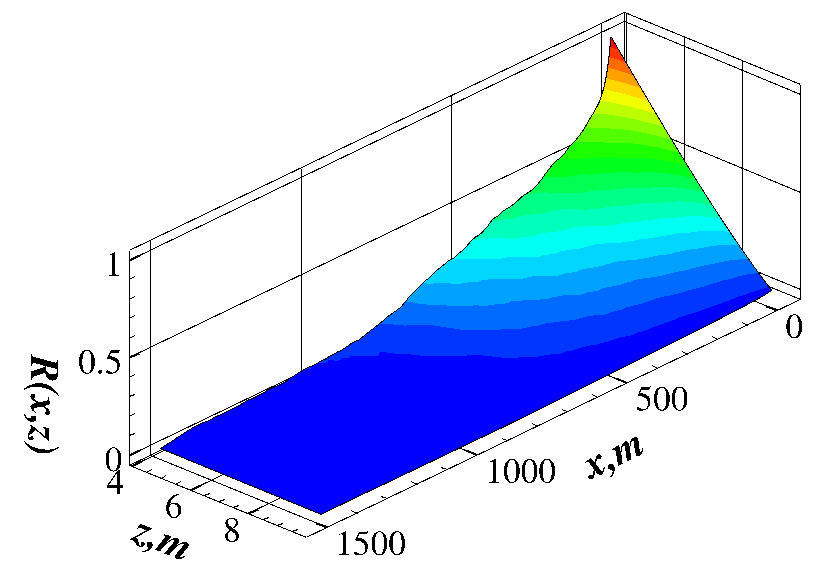} & 
    \includegraphics[width=0.3\textwidth]{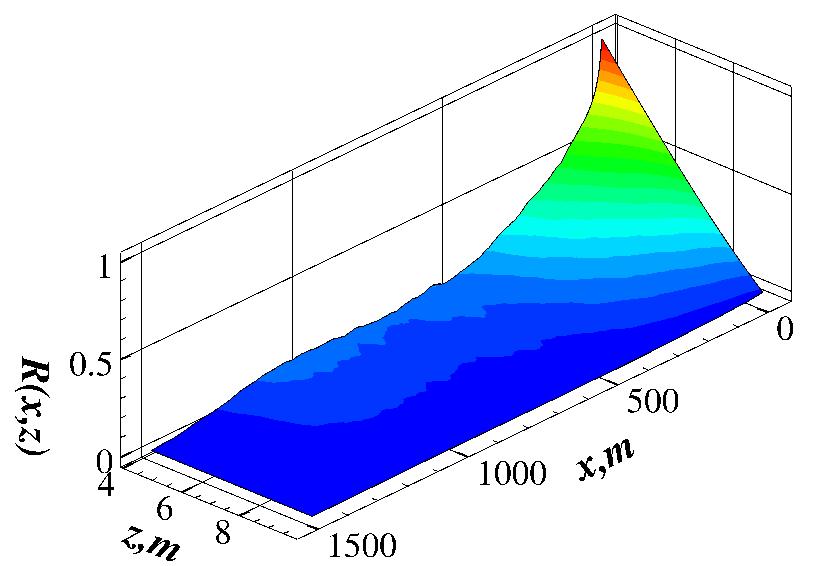} \\
    \includegraphics[width=0.3\textwidth]{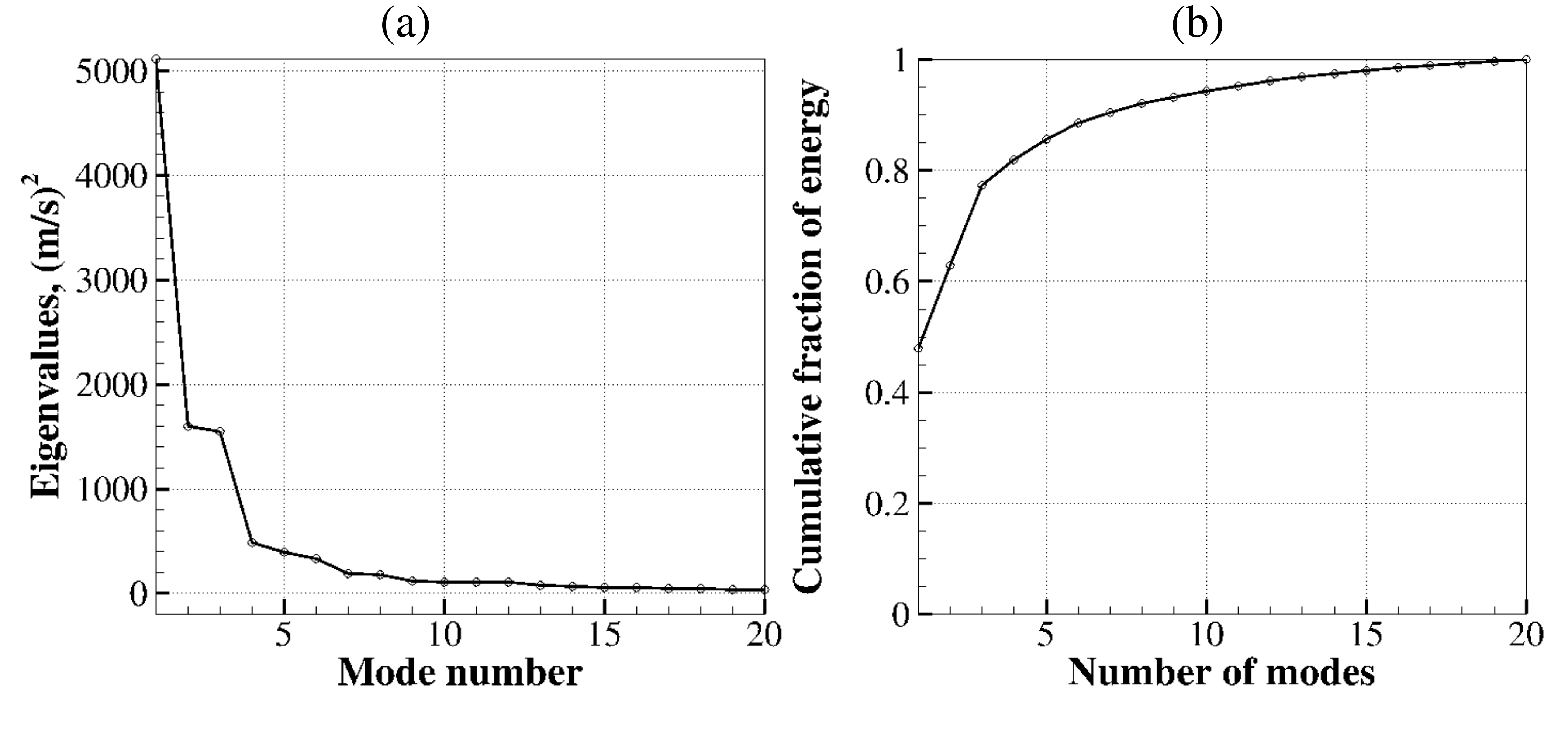} &
    \includegraphics[width=0.3\textwidth]{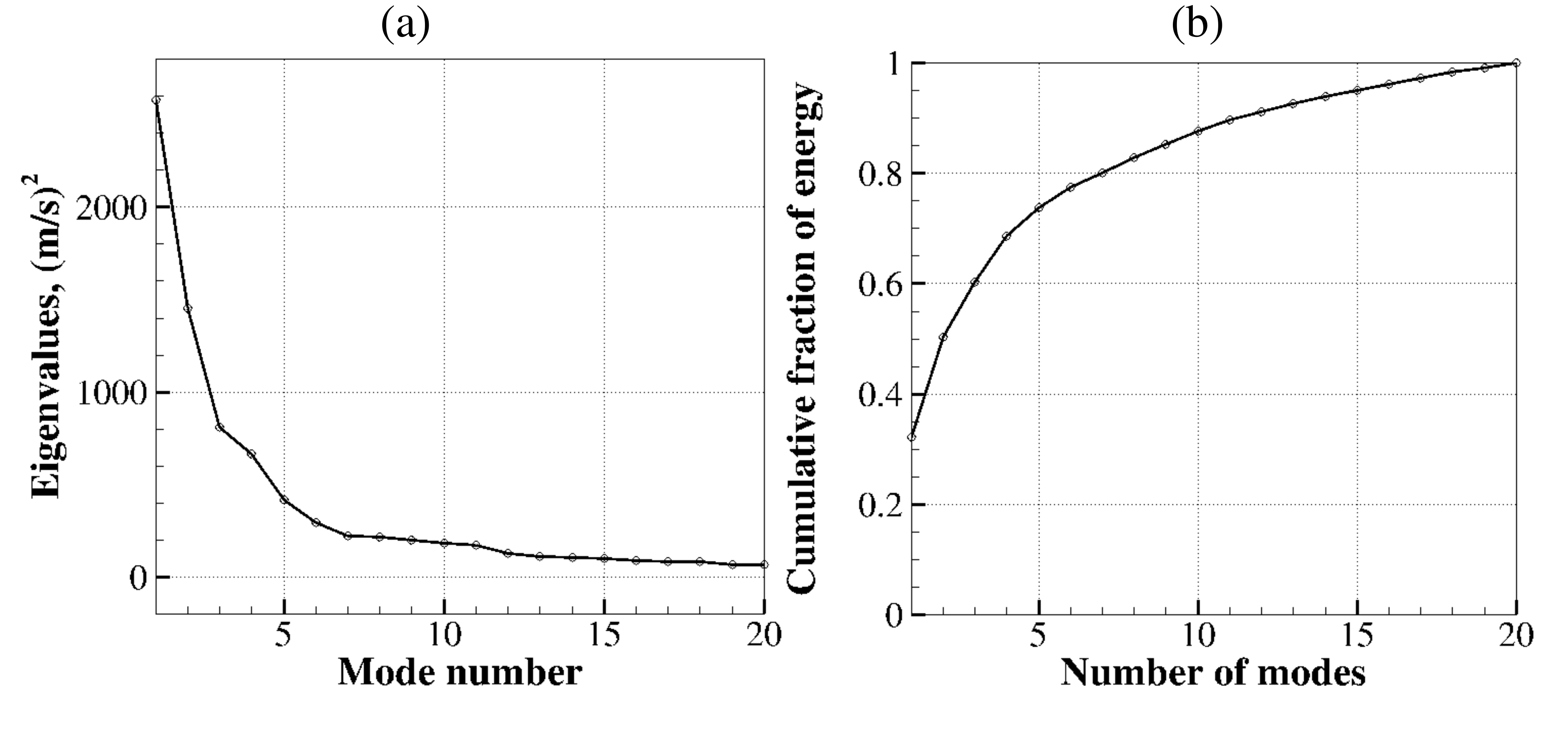} &
    \includegraphics[width=0.3\textwidth]{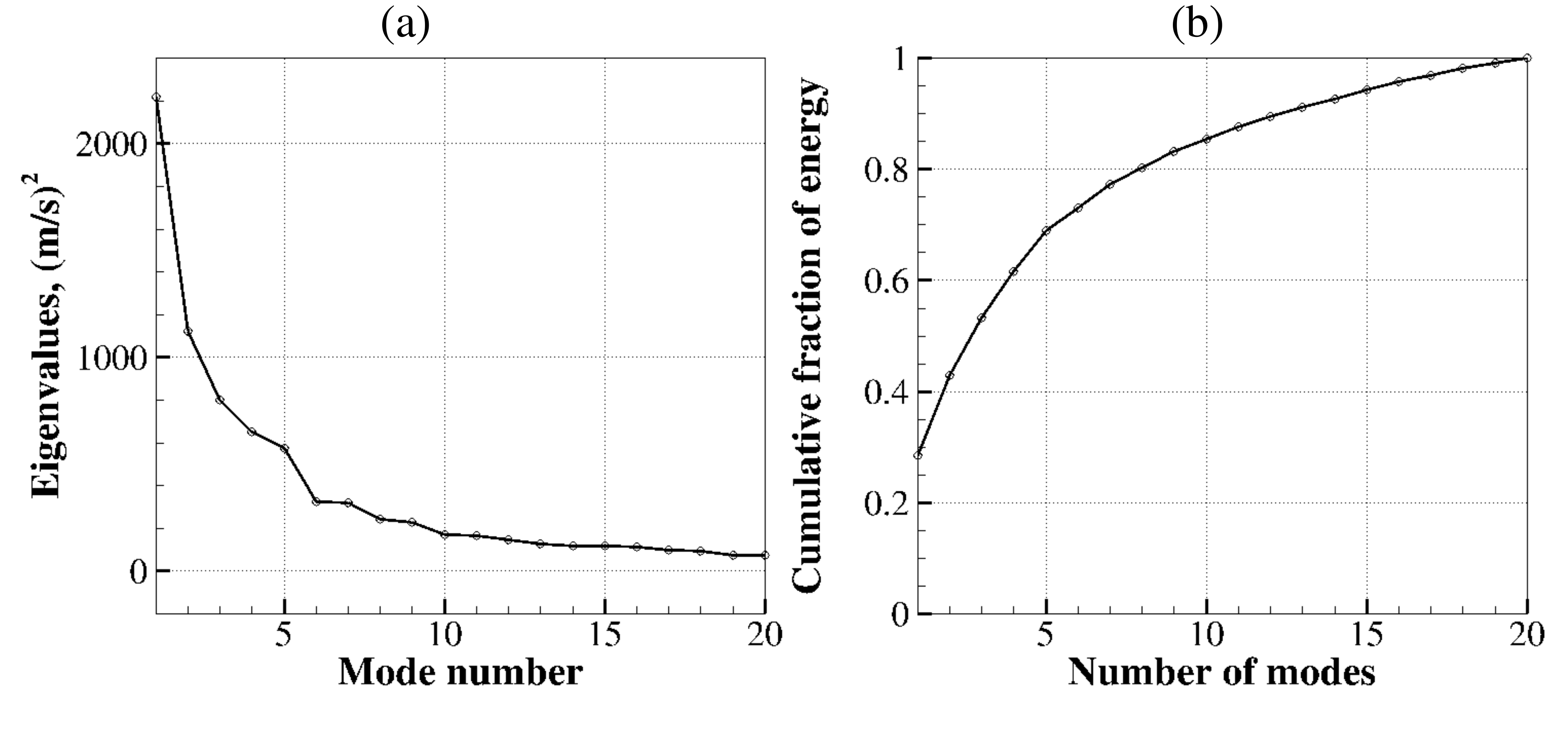} \\
    \includegraphics[width=0.3\textwidth]{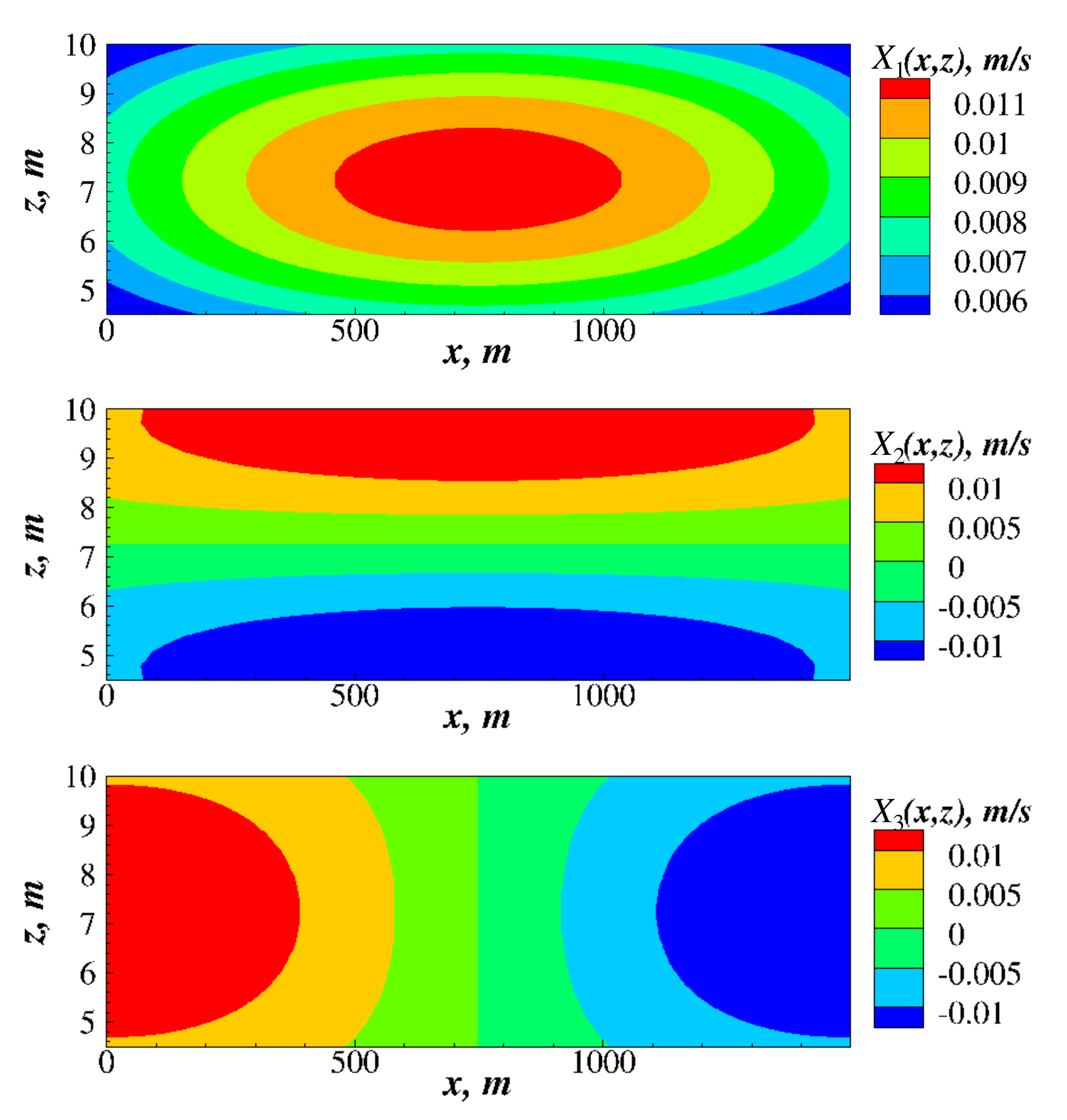} &
    \includegraphics[width=0.3\textwidth]{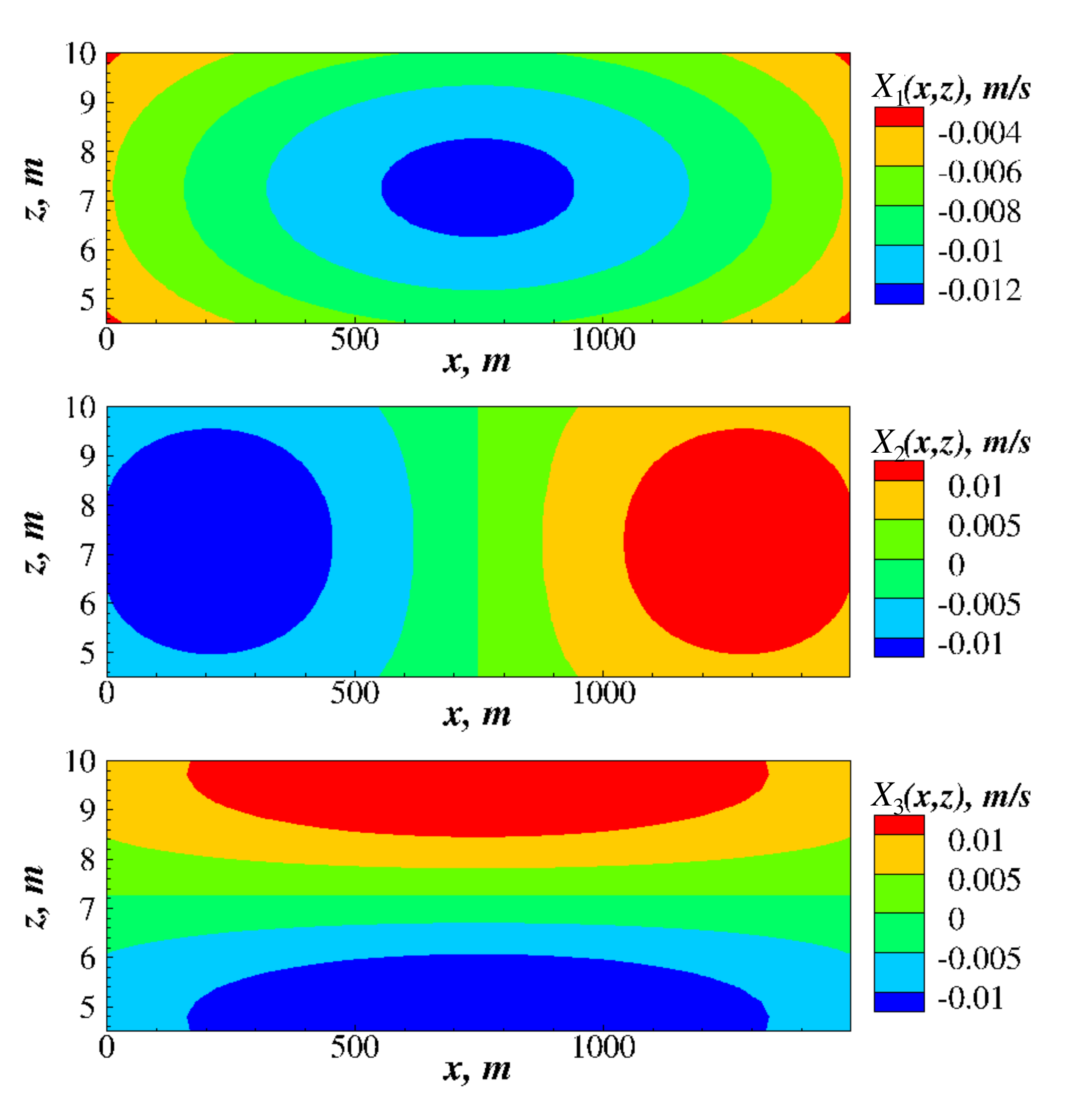} &
    \includegraphics[width=0.3\textwidth]{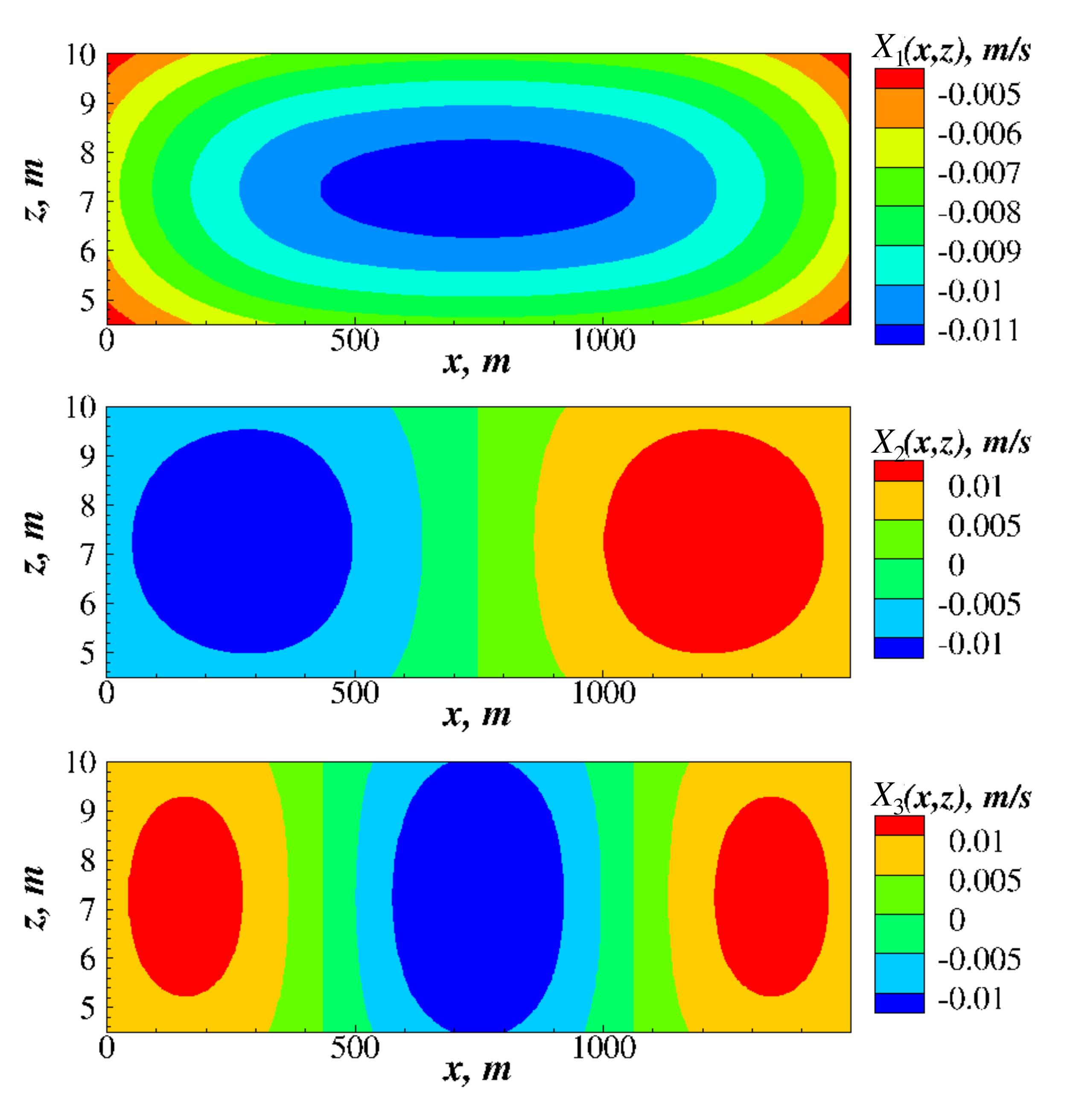} \\
    \includegraphics[width=0.3\textwidth]{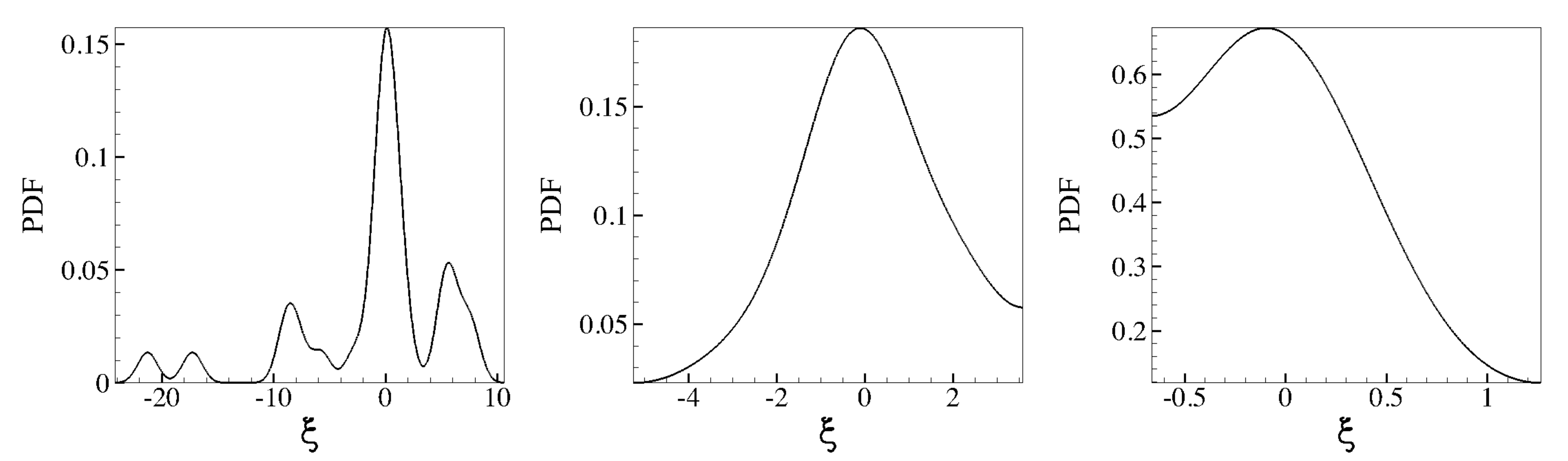} &
    \includegraphics[width=0.3\textwidth]{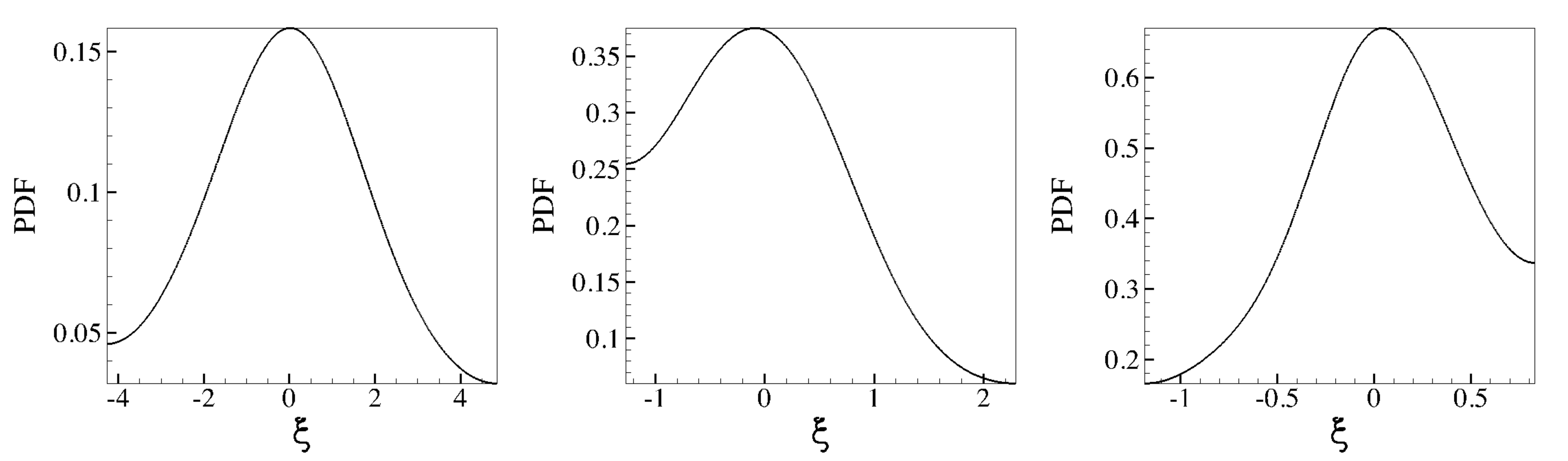} &
    \includegraphics[width=0.3\textwidth]{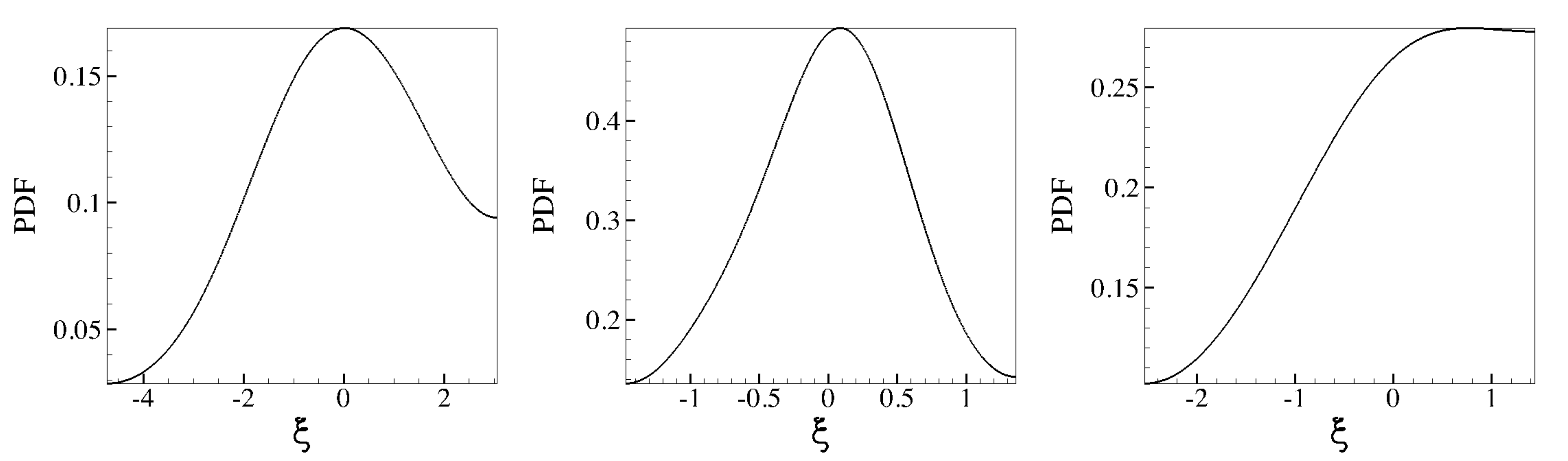} \\
    \hline
  \end{tabular}
  \caption{KLE results of the first three spatial modes} \label{tab:KLE}
\end{sidewaystable}

\subsection{The low-complexity wind model}\label{low-complexity model}
\begin{figure}[htbp]
\centering
 \includegraphics[width=0.9\textwidth]{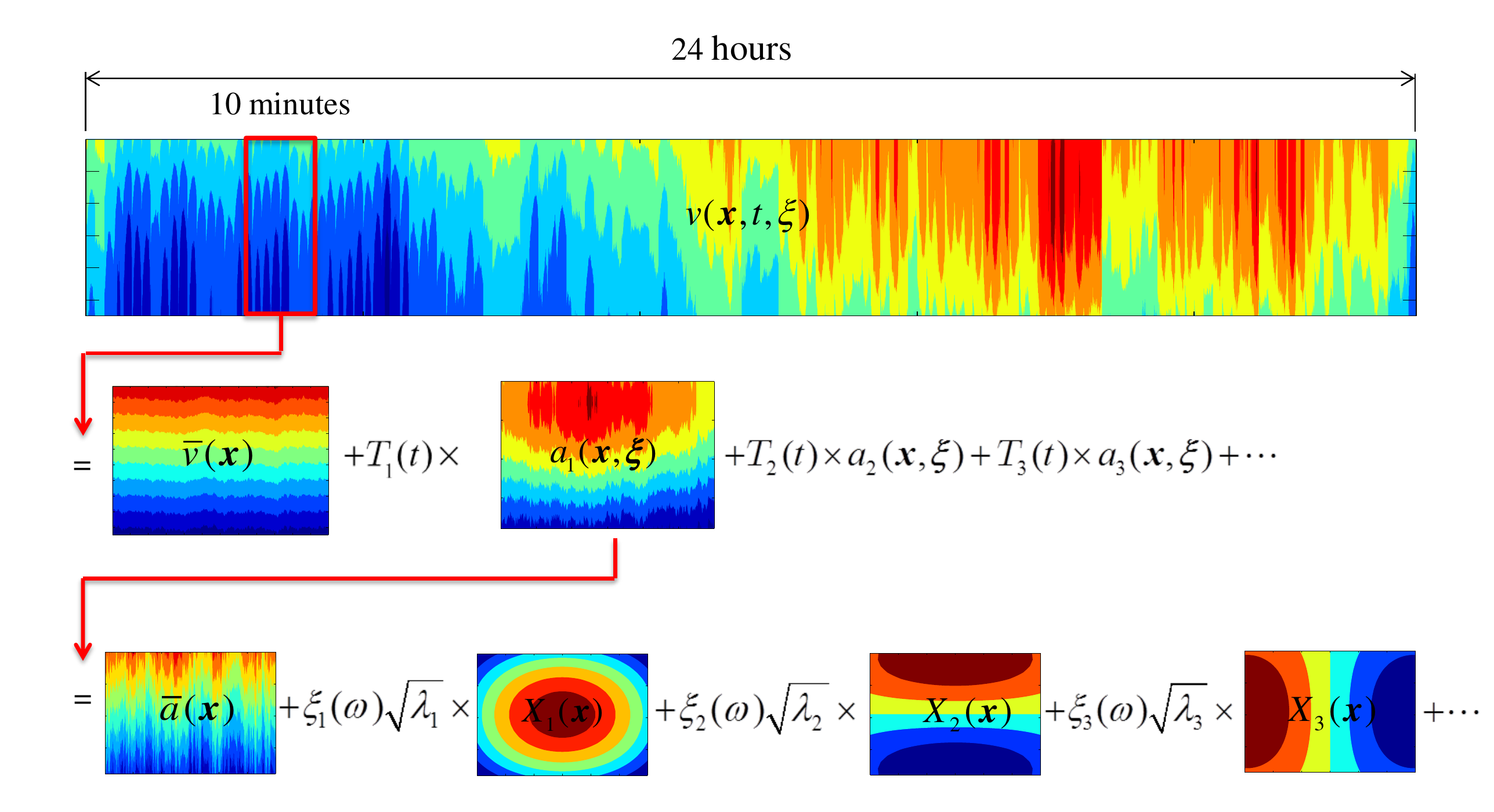}
\caption{The low-complexity model: Process of constructing synthetic wind snapshots}
\label{constsyn}
\end{figure}

The above two stages result in the calculation of the temporal functions, $T_i(t)$ (from Bi-orthogonal Decomposition), the spatial functions ($\bm X^i_j(\bm x)$) (from KLE decomposition), and the uncorrelated random variables, $\xi^i_j$ (from the KDE fitting).\footnote{The superscript in $\bm X^i_j(\bm x)$ and $\xi^i_j$ represents the index of Bi-orthogonal Decomposition terms whereas the subscript denotes the index of KLE terms.} Putting it all together gives us the low-complexity model for the wind
\begin{equation}
\bm u(\bm x,t,\bm\xi) =\sum K_{ij}\ T_i(t) \bm X^i_j(\bm x) \xi^i_{j}
\end{equation}
Note that $\bm X^i_j$ and $T_i$ are deterministic functions that encode spatial and temporal correlations of the wind. Different realizations (or stochasticity) of the wind is included into the low-complexity model via the uncorrelated random variables, $\xi^i_j$. Note that the probability distributions of $\xi^i_j$ are constructed in a data-driven way from the meteorology data. As we will show in the results section, only a few terms ($M=3$) are required to reconstruct a synthetic wind snapshot that \textbf{contains all the temporal and spatial correlations exhibited by the original data}. Fig.~\ref{constsyn} gives a graphical description of this process. Synthetic wind snapshots exactly mimicking the meteorological data can be constructed by sampling from the distributions of the random variables, $\xi^i_j$. Interestingly, only 800 KB of storage space is needed to store all the necessary information of the reduced-order model, whereas more than 200 MB of storage space is needed to store all the wind flow 
snapshots that are used in the analysis! This demonstrates the advantage of the low-complexity model in terms of data size. 

\subsection{Statistical comparison}\label{statistical comp}
As discussed previously, to get 90\% representation accuracy, five terms in Bi-orthogonal Decomposition and seven terms in Karhunen-Lo\`{e}ve Expansion are needed. However, for the purpose of demonstration, 1-, 3-, and 10-term BD and three terms in KLE for each BD term were used in the analysis. In order to quantify the accuracy, 28 realizations (same as the number of samples of meteorological data set) of the synthetic wind flow are generated. Each realization is a 24-hour wind data consisting of 144 ten-minute snapshots (see Fig.~\ref{constsyn} for one 24 hour synthetic dataset). 

\begin{figure}[h]
 \centering
  \includegraphics[width=0.6\textwidth]{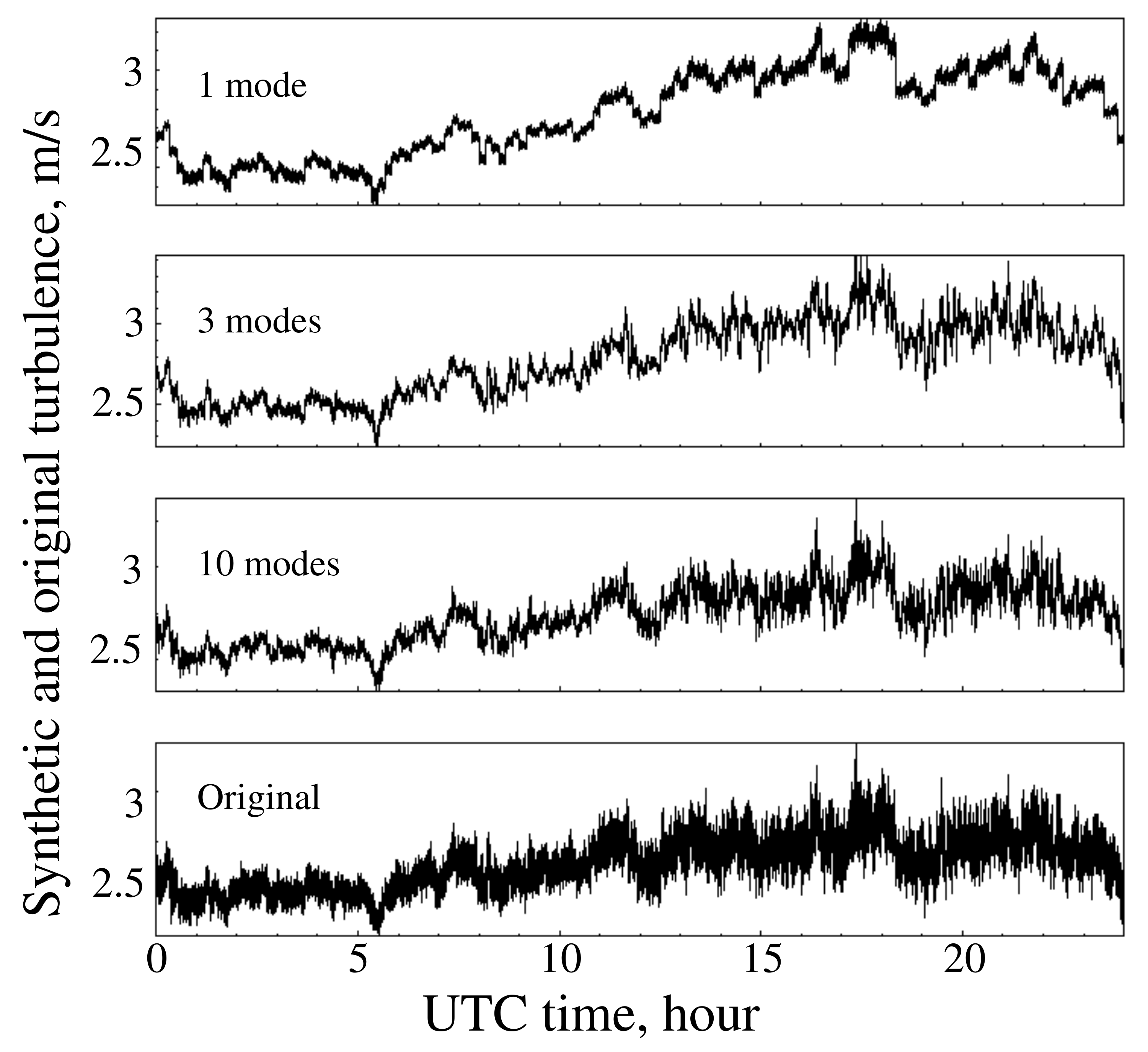}
 \caption{realizations of the stochastic flow}
 \label{realizations}
\end{figure}

In Fig.~\ref{realizations}, realizations of the stochastic flow at height 10m using one, three, and ten spatial modes are compared. Result shows that when smaller number of terms are used, the low frequency characteristics of the turbulence can be represented accurately. However, in order to describe high frequency behavior, higher modes in BD must be used. In the third plot that uses ten temporal modes, diurnal fluctuation in the variance of turbulence is addressed by including the second and third temporal modes (see Fig.~\ref{eigenfunc}).

\begin{figure}[h]
 \centering
   \includegraphics[width=0.95\textwidth]{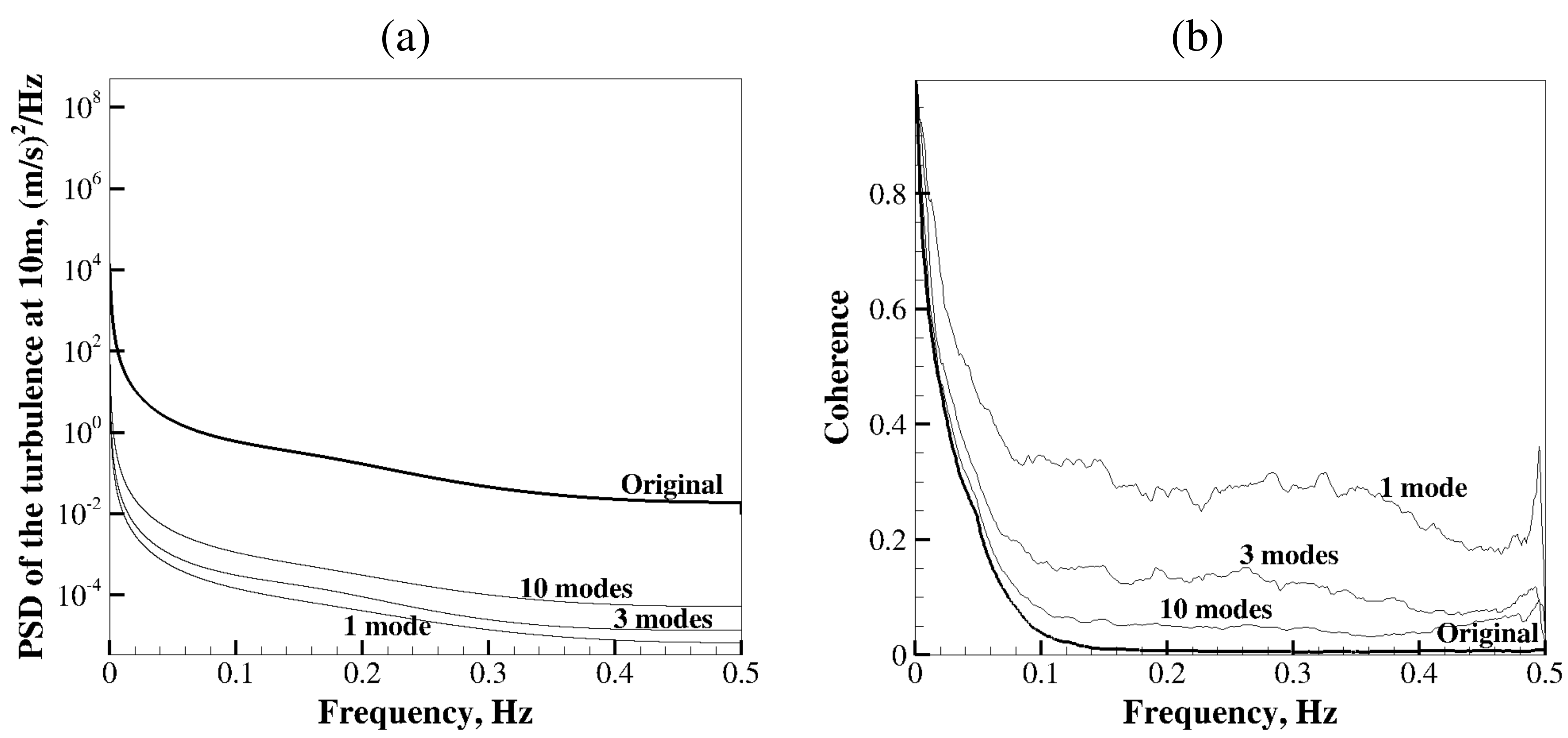}
   \caption{PSDs and coherences of the synthetic turbulence at 10 m using one, three, and ten modes compared with original flow}
   \label{psdncoh}
\end{figure}

One of the most important variables is the power spectral density (PSD) that describes how the turbulent energy is distributed across the frequency spectrum. To verify the similarity of the reconstructed wind flow and the meteorological data, their PSD functions are compared in Fig.~\ref{psdncoh} (a). The figure shows that the synthetic wind flow accurately reproduces energy in low frequency region. However, more terms need to be used in the representation to preserve energy in high frequency region. It is worth mentioning that even with only one BD mode, the synthetic wind mimics the true data in great detail.

Coherence spectrum is another important variable that describes the similarity of turbulence at two different spatial locations. Comparison of coherences between same set of time-series used in PSD analysis is given in Fig.~\ref{psdncoh} (b). Note that even by using only one BD mode, the synthetic wind can still preserve most of the spatial coherence information. This can be explained by the advantages of STD. Since STD performs orthogonal decompositions of both temporal and spatial covariances, the temporal variations and spatial correlations are preserved in an optimal way. 

\begin{figure}[htbp]
 \centering
  \includegraphics[width=0.95\textwidth]{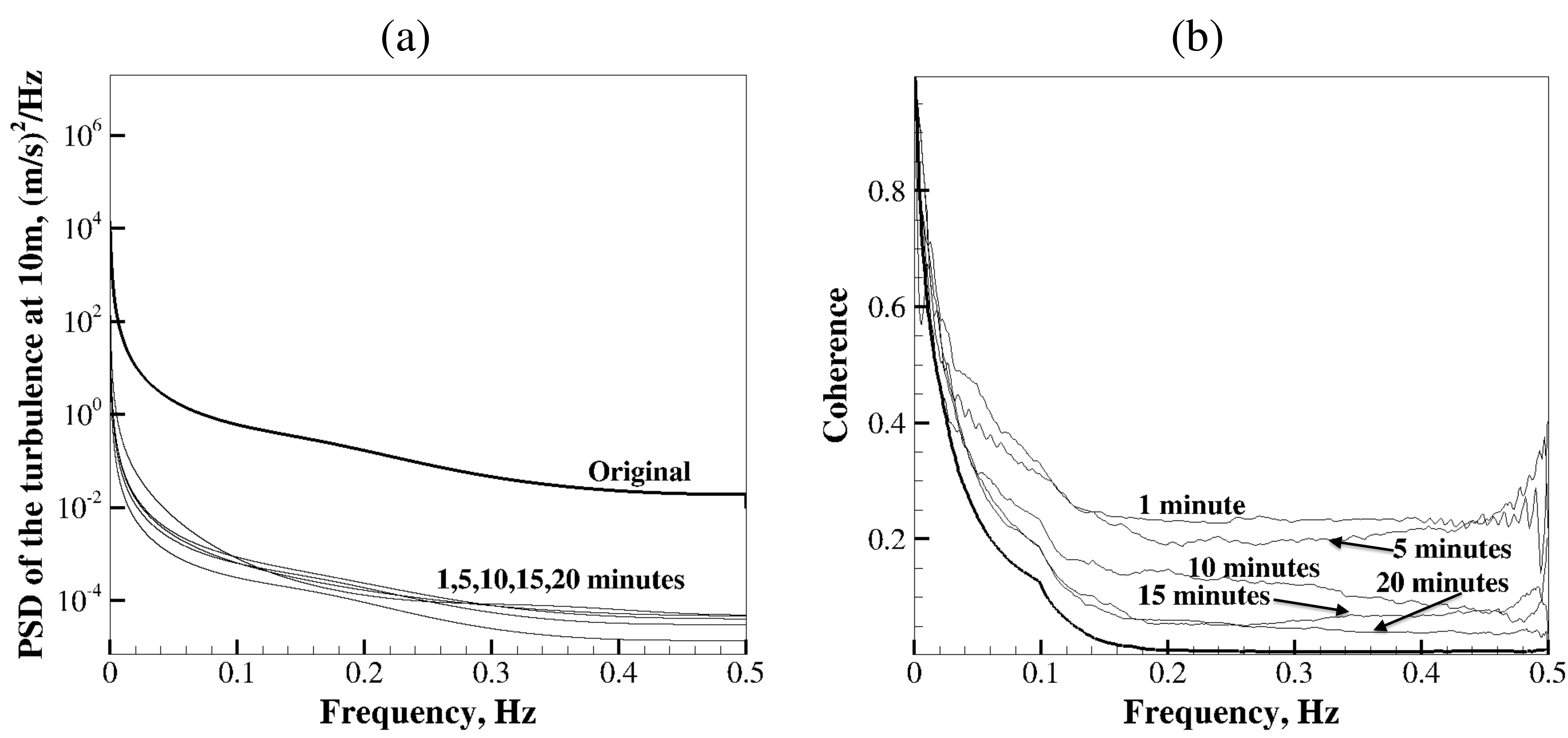}
 \caption{Comparing PSDs and coherences of the synthetic turbulence using different choices of snapshot} \label{interval choices}
\end{figure}

As previously mentioned, there are different choices of time intervals when constructing wind flow snapshots. In above content, only analysis on 10 minute snapshot was performed. The reason for choosing 10 minute snapshot can be demonstrated by comparing PSDs and coherences of the synthetic turbulence using different choices of time interval averaging processes. Fig.~\ref{interval choices} shows results of such analysis. Since the accuracy of PSDs of synthetic wind only depends on how many total terms in the decomposition, there is no apparent difference in PSDs of different choices of time intervals. On the other hand, comparison of coherences tells us if performing analysis on snapshots that are constructed from less than 10 minute field, important information in turbulence coherence will be lost. Coherences correspond to 10, 15, and 20 minute snapshots are very close, which means the additional turbulence data provided by 15 and 20 minute snapshots, comparing to 10 minute snapshot, does not contain more 
coherence information. This result is consistent with GL guideline that states the coherence structure in turbulence can be assumed to be unchanged in 10 minutes period of time~\cite{GermanischerLloyd2010}.

As illustrated in Fig.~\ref{eigenvalue&sm}, a certain number of spatial modes are needed to achieve a desired representation accuracy. On the other hand, the number of random inputs to the reduced order model equals the number of spatial points on each snapshot. By using different choices of time intervals, the resulting snapshots have different dimensions. Table~\ref{number relationship} shows the relationship between the number of random inputs and the number of needed spatial modes for certain level of accuracy. Based on the table, when the number of random inputs increases by 20 times, the needed random spatial modes to achieve 94\% of accuracy only increases from 3 to 8. In other words, the required number of spatial modes does not strongly depend on the number of random inputs, which suggest that this approach (TSD) is practical for problems with very large stochastic dimensions.
\begin{table}[htbp]%
\caption{Relationship between the number of random spatial modes (94\% accuracy) and the number of random inputs.}
\centering
\begin{tabular}{cccccc}
\hline %
Number of random inputs&1200&6000&12000&18000&24000\\
Number of spatial modes&3&5&6&7&8\\\hline
\end{tabular}
\label{number relationship}
\end{table}

The temporal covariance function of synthetic data $\hat{C}(t,t')$ is constructed. Comparing it with the temporal covariance of the original data reveals that they have almost identical pattern. 
\begin{figure}[htbp]
\centering
 \includegraphics[width=0.95\textwidth]{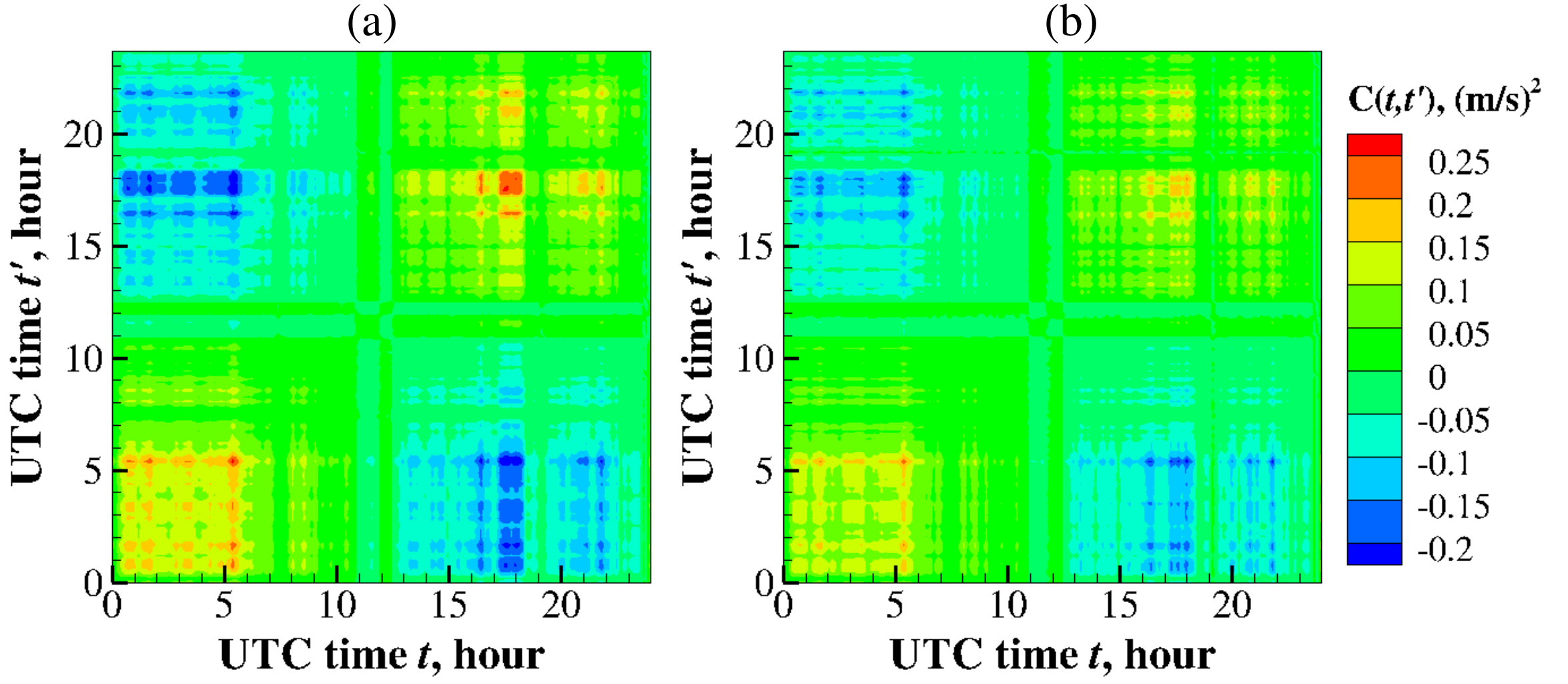}
\caption{Comparison of covariance functions of original (a) and synthetic (b) flows}
\label{syncor}
\end{figure}
The temporal error (or the information loss) can be defined as the $L^2$-norm of the difference between the covariance functions: 
\begin{equation}
 \epsilon = \frac{\parallel C-\hat{C} \parallel_2}{\parallel C \parallel_2} = \left( \frac{\sum_{i=1}^n \sum_{j=1}^n \mid C(t_i,t_j)-\hat{C}(t_i,t_j) \mid^2}{\sum_{i=1}^n \sum_{j=1}^n \mid C(t_i,t_j) \mid^2} \right)^{1/2}.
\end{equation}
where $n=144$ is the number of snapshots in 24 hours. The temporal information-loss using a 9-term expansion is 2.25\%. Thus, a 9-term data-driven expansion reproduces the temporal and spatial covariance of the original meteorological data to 97.75\% accuracy.

\section{Discussion and Conclusion}\label{conclusion}
Incorporating the effects of randomness in wind is critical for a variety of application involving wind energy. Continuous advances in data-sensing and meteorology has made possible the availability of large data-sets of location-, topography-, diurnal-, and seasonal sensitive meteorology data. While this data contains rich information, ease of use is bottlenecked by the unwieldy data-sizes. A pressing challenge is to utilize this data to construct a location-, topography-, diurnal-, and seasonal-, dependent low-complexity model that is easy to use and store. We formulate a data-driven mathematical framework that is capable of representing the spatial- and temporal- correlations as well as the inherent randomness of wind into a low-complexity parametrization. We leverage data-driven decomposition strategies like Bi-orthogonal Decomposition and Karhunen-Lo\`{e}ve expansion for constructing the low-complexity model. We provide a software package that can be used to construct the low-complexity model to the community. 

There are few points should be emphasized. First, the meteorological data used in this analysis is measured by only one met tower. Because the lack of the measurement on transverse direction, getting wind filed snapshots on the plane of turbine rotor becomes impossible. To circumvent this insufficiency in measurement, wind filed snapshots on the plane that is perpendicular to the rotor is used, for which certain time interval has to be chose to construct snapshots. It is worth noting that the framework is generally applicable to a variety of meteorology data, and the its applicability should not be affected by choosing different snapshot constructing time intervals. In addition, the framework is able to incorporate both short-term (10-minute) and long-term (years) temporal coherences as long as corresponding data is available. Third, while the mathematical framework developed here is used to analyze wind speed, it can also be used to represent other atmospheric data such as temperature and carbon dioxide flux. This framework can also be naturally extended to represent ocean waves, which is crucial for off-shore wind turbine siting, layout and design analysis.

\section{Acknowledgement}
BG and QG gratefully acknowledge NSF 1149365. 

\section{References}
\bibliographystyle{model1-num-names}
\bibliography{mybib}


\section{Appendix}
\begin{appendices}
\numberwithin{equation}{section}

\section{Derivation of Minimizing Representation Error using Euler Lagrange Equation}\label{reperror}
\newtheorem{theorem}{Theorem}
\begin{theorem}
\emph{(Euler's equation)}
\label{Euler}
 Let $J[y]$ be a functional of the form
  \begin{equation}
  \int_a^b F(x,y,y')dx,
  \end{equation}
 defined on the set of functions $y(x)$ which have continuous first derivatives in $[a,b]$ and satisfy the boundary conditions $y(a)=A$, $y(b)=B$. Then a necessary condition for $J[y]$ to have an extremum for a given function $y(x)$ is that $y(x)$ satisfy Euler's equation
  \begin{equation}
  F_y-\frac{d}{dx} F_{y'}=0 
  \end{equation}
\hfill \qed
\end{theorem} 

The goal of BD is to minimize the error functional
\begin{equation}
 \mathscr{E}[T_1,\cdots,T_M] = \int_T \langle \varepsilon, \varepsilon \rangle_X dt.
\end{equation}
Substituting representation error (Eqn.~\ref{error}) in above equation yields
\begin{equation}
\begin{split}
 & \mathscr{E}[T_1,\cdots,T_M] \\
 & = \int_T \langle u(\bm x,t,\bm\xi) - \sum_{i=1}^M K_i \bm \varPhi_i(\bm x,\bm\xi)T_i(t),\ u(\bm x,t,\bm\xi) - \sum_{j=1}^M K_i \bm \varPhi_j(\bm x,\bm\xi)T_j(t) \rangle_X\ dt 
\end{split}
\end{equation}
According to the definition of spatial inner product (Eqn.~\ref{spatial inner}), the error functional becomes

\begin{equation}
\begin{split} 
 & \mathscr{E}[T_1,\cdots,T_M] \\
 & = \int_T \int_{\bm X} \left[ u(\bm x,t) - \sum_{i=1}^M K_i\bm \varPhi_i(\bm x)T_i(t) \right]\left[\ u(\bm x,t) - \sum_{j=1}^M K_j\bm \varPhi_j(\bm x)T_j(t)\right]\ d\bm x\ dt \\
 & = \int_T \int_{\bm X} u^2(\bm x,t)-2u(\bm x,t)\sum_{i=1}^M K_i\bm \varPhi_i(\bm x)T_i(t) \\
 & \quad +\sum_{i=1}^M K_i\bm \varPhi_i(\bm x)T_i(t)\sum_{j=1}^M K_j\bm \varPhi_j(\bm x)T_j(t)\ d\bm x\ dt \\
 & = \int_T \int_{\bm X} u^2(\bm x,t)-2u(\bm x,t)\sum_{i=1}^M K_i\bm \varPhi_i(\bm x)T_i(t) +\sum_{i=1}^M K^2_i\bm \varPhi_i^2(\bm x)T_i^2(t)\ d\bm x\ dt \\
 & = \int_T\left[\int_{\bm X} u^2(\bm x,t)d\bm x-2\sum_{i=1}^M K_i T_i(t)\int_{\bm X} u(\bm x,t)\bm \varPhi_i(\bm x)d\bm x +\sum_{i=1}^M K^2_i T_i^2(t)\right]\ dt \\
 & = \int_T\ F(t, T_1,\cdots,T_M, T'_1,\cdots,T'_M) dt
\end{split}
\end{equation}
where
\begin{equation}
\begin{split}
 & F(t, T_1,\cdots,T_M, T'_1,\cdots,T'_M) \\
 & = \int_{\bm X} u^2(\bm x,t)d\bm x-2\sum_{i=1}^M K_i T_i(t)\int_{\bm X} u(\bm x,t)\bm \varPhi_i(\bm x)d\bm x+\sum_{i=1}^M K_i^2 T_i^2(t)
\end{split}
\end{equation}
In above derivation, the orthogonality of basis functions $\bm \varPhi_i$ and $T_i$ is applied. According to Theorem~\ref{Euler}, the error functional $\mathscr{E}[T_1,\cdots,T_M]$ has extremum when $T_i$ satisfy Euler's equations
\begin{equation}
 F_{T_i}-\frac{d}{d\bm x}F_{T'_i} = 0.
\end{equation}
That is
\begin{equation}
 \int_{\bm X}u(\bm x,t)\bm \varPhi_i(\bm x)d\bm x-K_i T_i(t)=0,
\end{equation}
which can be simplified as
\begin{equation}
 T_i(t)=\frac{1}{K_i}\langle u(\bm x,t,\bm\xi),\bm \varPhi_i(\bm x,\bm\xi) \rangle_X. \label{bestT}
\end{equation}

Applying temporal inner product $\langle \cdot, T_i \rangle_T$ to Bi-orthogonal Decomposition (Eqn.~\ref{BD}) and considering the orthogonality of the temporal basis functions yields
\begin{equation}
 \bm \varPhi_i(\bm x,\bm\xi) = \frac{1}{K_i} \langle u(\bm x,t,\bm\xi), T_i(t) \rangle_T. \label{bestPhi}
\end{equation}
Above two equations define a coupled relationship between $T_i$ and $\bm \varPhi_i$. Now that we have two unknowns and two equations, they can be solved by substituting Eqn.~\ref{bestPhi} into Eqn.~\ref{bestT}, which results in eigenvalue problem for temporal modes
\begin{equation}
 K^2_i T_i(t) = \int_T C(t,t') T_i(t')dt',  \label{BD eigen problem 1}
\end{equation}
where $C(t,t')$ is called temporal covariance that can be obtained by taking the inner product in spatial domain, i.e
\begin{equation}
 C(t,t') = \langle u(\bm x,t,\bm\xi),u(\bm x,t',\bm\xi) \rangle_X.
\end{equation}
Setting $\mu_i=K^2_i$, Eqn.~\ref{BD eigen problem 1} can be transformed to
\begin{equation}
 \mu_i T_i(t) = \int_T C(t,t') T_i(t')dt',  \label{BD eigen problem 2}
\end{equation}
where $\mu_i$ and $T_i(t)$ are eigenvalues and eigenfunctions of the covariance function $C(t,t')$. The optimal choice of temporal modes $T_i$ and $\bm \varPhi_i$ can be obtained by solving the eigenvalue problem.

\section{Numerical solution to the Generalized Eigenvalue Problem}\label{numericaleigen}
We want to find the numerical solution of eigen-equation:
\begin{equation}
 \int_{\bm{X}} C(\bm x_1,\bm x_2) \, f_{i}(\bm x_1) \, d\bm x_1 = \lambda_{i} \, f_{i}(\bm x_2).
\end{equation}
To this end, eigenfunction can be approximated by linear combination of $N$ basis functions
\begin{equation}
 f_{k}(\bm x) = \sum^N_{i=1} \, d^{(k)}_{i} \, h_{i}(\bm x).
\end{equation}
Substitute above equation to the eigen-equation and set the error to be orthogonal to each basis function yields
\begin{equation}
 \sum^N_{i=1} d^{(k)}_{i} \left[ \int_{\bm{X}} \left[ \int_{\bm{X}} C(\bm x_1,\bm x_2) h_{i}(\bm x_2) d\bm x_2 \right] h_{j}(\bm x_1) d\bm x_1 - \lambda_{n} \int_{\bm{X}} h_i(\bm x)h_j(\bm x) d\bm x \right] = 0.
\end{equation}
Above equation can be written in matrices form
\begin{equation}
 A D = B D \Lambda.
\end{equation}
\begin{equation}
 A_{ij} = \int_{\bm{X}} \int_{\bm{X}} C(\bm x_1,\bm x_2) h_{i}(\bm x_1) h_{j}(\bm x_2) d\bm x_1 d\bm x_2.
\end{equation}
\begin{equation}
 B_{ij} = \int_{\bm{X}} h_i(\bm x)h_j(\bm x) d\bm x = \int_{\bm{X}} H^T(\bm x) H(\bm x) d\bm x.
\end{equation}
\begin{equation}
 D_{ij} = d^{(j)}_i.
\end{equation}
\begin{equation}
 \Lambda_{ij} = \delta_{ij} \lambda_i.
\end{equation}
where  $H(\bm x) = \left( h_1(\bm x), h_2(\bm x), \cdots, h_N(\bm x) \right)$. Matrix A can be rewritten as
\begin{align*} 
A & = \int_{\bm{X}} \int_{\bm{X}} H^T(\bm x_1) C(\bm x_1,\bm x_2) H(\bm x_2) d\bm x_1 d\bm x_2 \\
& =\int_{\bm{X}} \int_{\bm{X}} H^T(\bm x_1) H(\bm x_1) C H^T(\bm x_2) H(\bm x_2) d\bm x_1 d\bm x_2 \\
& =\int_{\bm{X}} H^T(\bm x_1) H(\bm x_1) d\bm x_1 \, C \, \int_{\bm{X}} H^T(\bm x_2) H(\bm x_2) d\bm x_2 \\
& =B C B.
\end{align*}
and
\begin{align*} 
C(\bm x_k,\bm x_l) & = \sum^N_{i=1} \sum^N_{j=1} h_i(\bm x_k) C_{ij} h_i(\bm x_l) \\
& = h_k(\bm x_k) C_{kl} h_l(\bm x_l) \\
& = C_{kl}.
\end{align*}

\end{appendices}

\end{document}